%% file: Collas,_Dettweiler,_Reiter,_Sawin_-_Elliptic_Convolution_G2-motives.tex
\documentclass[a4paper,11pt, english,twoside,centertags,leqno]{article}

	\usepackage[utf8]{inputenc}
	\usepackage[T1]{fontenc}

	\usepackage{amsmath, amssymb,latexsym,amsxtra,amscd,amsthm}
	\usepackage{mathtools}
	\usepackage{stmaryrd} 
	\usepackage[all]{xy}
	\usepackage{rotating}

	\usepackage[lighttt]{lmodern} 
	\usepackage[protrusion=true,expansion=true,kerning, spacing]{microtype}
	\usepackage{listings} 
	\usepackage{graphicx}
	\usepackage{calc} 
	\graphicspath{{Figures/}}
	\usepackage{subcaption}
	\usepackage[usenames,dvipsnames,svgnames]{xcolor}
	\usepackage[shortlabels]{enumitem} 

	\usepackage[english]{babel}
	\usepackage[textwidth=15cm,textheight=21cm,heightrounded,hcentering]{geometry}
	\usepackage{titlesec}
	\titleformat{\section}[block]{\large\bfseries}{\thesection.}{0.5em}{}	
	\titleformat{\subsection}[runin]{\normalfont\bfseries}{\thesubsection.}{0.5em}{}[.]

	\usepackage{caption}
	\usepackage{ctable} 
	\usepackage{multirow}
	
	\theoremstyle{plain}
	        \newtheorem{thm}{Theorem}[subsection]
	        \newtheorem{cor}[thm]{Corollary}
	        \newtheorem{lem}[thm]{Lemma}
	        \newtheorem{prop}[thm]{Proposition}

	\theoremstyle{definition}	
	        \newtheorem{defn}[thm]{Definition}

	\theoremstyle{remark}
	        \newtheorem{rem}[thm]{Remark}

	\numberwithin{equation}{section}
	
	\def\AA{\mathbb{A}}
	\def\CC{\mathbb{C}}
	\def\GG{\mathbb{G}}
	\def\NN{\mathbb{N}}
	\def\pr{\mathbb{P}}
	\def\QQ{\mathbb{Q}}
	\def\RR{\mathbb{R}}

	\newcommand{\x}{{\bf x}}

	\let\oldvee\vee
	\renewcommand{\vee}{{\scriptscriptstyle\oldvee}}

	\let\epsilon\varepsilon
	\def\loccit{loc.\kern3pt cit.{}\xspace}
	\def\cf{cf.\kern.3em}
	\def\eg{e.g.\kern.3em}
	
	\def\resp{\text{resp.}\kern.3em}
	\let\setminus\smallsetminus
	\let\leq\leqslant
	\let\geq\geqslant

	\let\tilde\widetilde
	\let\bar\overline
	\let\hat\widehat

	\newcommand{\ptbl}{.\kern .2em }
	\newcommand{\an}{\mathrm{an}}

	\DeclareMathOperator{\im}{Im}

	\DeclareMathOperator{\GL}{GL}
	\DeclareMathOperator{\SO}{SO}

	\DeclareMathOperator{\Hom}{Hom}
	\DeclareMathOperator{\rk}{rk}

	\newcommand{\SL}{{\rm SL}}

	\renewcommand{\pr}{{\frak p}}

	\newcommand{\To}{\;\longrightarrow\;}

	\renewcommand{\pr}{{{\rm pr}}}

	\newcommand{\Perv}{{\rm Perv}}

	\renewcommand{\char}{{\rm char}}

	\newcommand{\PP}{{\mathbb{P}}}

	\newcommand{\y}{{\bf{y} }}

	\newcommand{\bQl}{\bar{\mathbb{Q}}_\ell}
	
	\def\to{\mathchoice{\longrightarrow}{\rightarrow}{\rightarrow}{\rightarrow}}

	\def\To#1{\mathchoice{\xrightarrow{\textstyle\kern4pt#1\kern3pt}}{\stackrel{#1}{\longrightarrow}}{}{}}

	\def\isoTo#1{\xrightarrow[\sim]{\textstyle\kern4pt#1\kern3pt}}

	\newcommand{\RedefinitSymbole}[1]{%
	\expandafter\let\csname old\string#1\endcsname=#1
	\let#1=\relax
	\newcommand{#1}{\csname old\string#1\endcsname\,}%
	}
	\RedefinitSymbole{\forall} \RedefinitSymbole{\exists}
	
	\usepackage{hyperref} 

	\title{Monodromy of Elliptic Curve Convolution, Seven-point Sheaves of $G_2$-type, and Motives of Beauville Type}

	\author{Benjamin \textsc{Collas}\footnote{B.~Collas: 
	Department of Mathematics,
	University of Bayreuth,
	95440 Bayreuth,
	Germany;
	{benjamin.collas@uni-bayreuth.de} -- supported by DFG grant SPP1786}, 
	Michael \textsc{Dettweiler}\footnote{M.~Dettweiler: Department of Mathematics,
	University of Bayreuth,
	95440 Bayreuth,
	Germany;
	{michael.dettweiler@uni-bayreuth.de}},
	Stefan \textsc{Reiter}\footnote{S.~Reiter: Department of Mathematics,
	University of Bayreuth,
	95440 Bayreuth,
	Germany;
	{stefan.reiter@uni-bayreuth.de}}
	\ and Will \textsc{Sawin}\footnote{W.~Sawin: Department of Mathematics, Columbia University, 10027 New York, NY, USA; sawin@math.columbia.edu}
}

	\newcommand{\PPE}{\mathcal{P}_E}
	\newcounter{introPara}

	\date{\today}
\begin{document}

	\maketitle

\begin{abstract}
	We study the Tannakian properties of the category of perverse sheaves on elliptic curves endowed with the convolution product. We establish that for certain sheaves with unipotent local monodromy over seven points the corresponding Tannaka group is isomorphic to $G_2$. This monodromy approach generalizes a result of Katz on the existence of $G_2$-motives in the middle cohomology of deformations of Beauville surfaces.
\end{abstract}

\tableofcontents

\section*{Introduction}
\addcontentsline{toc}{section}{Introduction}
For an abelian variety $A,$
the {\it convolution} $K_1\ast K_2$ of two objects $K_1,K_2$ in $D^b_c(A)$ is given by the 
derived pushforward of $K_1\boxtimes K_2$ along the addition map $a\colon A\times A\to A$ by:
\[
K_1\ast K_2=Ra_*(K_1\boxtimes K_2).
\] 
As studied among others by Katz \cite{KatzEll}, Kr\"amer and Weissauer \cite{Weiss1}, \cite{KraemerWeissauer}, properties of the convolution on abelian varieties leads to a \emph{neutral Tannakian category} $(\mathcal{P}_A,\star)$ within the category $\Perv(A)$ of perverse sheaves over $A$ (see also Section \ref{secrev}); the category $\mathcal{P}_A$ contains all irreducible perverse sheaves of geometric origin, given by a \emph{variation of motives over $A$}, which are not shifts of  local systems on the whole of $A$. This paper considers the case of elliptic curves $E$ over an algebraic closed field $\bar k$ and deals with Tannakian properties of the convolution of specific \emph{seven-point sheaves} of $\PPE$. Original motivations come from a similar Tannakian constructions in genus $0$, i.e on $\GG_{\bar k}$, that lead to results by Deligne and Terasoma for \emph{periods of mixed Tate motives} \cite{DELRAC,TER06}, and to applications in \emph{Regular Inverse Galois Theory} \cite{dw}.

More precisely, on $\mathcal{P}_A$ the $H^0$-functor switches convolution and tensor products -- i.e. \[
H^0(A,K_1\ast K_2)=H^0(A,K_1)\otimes H^0(A,K_2),
\] 
and in the case where $L$ is an irreducible local system on a dense open subset $j\colon U\hookrightarrow A$ with $K=(j_{!*}L)[n]\in \mathcal{P}_A$ (where $dim(A)=n$), the fibre functor is the middle cohomology of $j_{!*}L$ on $A$: 
\[
H^0(A,K)\simeq H^{n}(A,j_{!*}L).
\] 

To each object $K\in \mathcal{P}_A$, one associates further by Tannakian formalism an affine group scheme $G_K$ which is a reductive subgroup of $\GL(H^0(A,K))$. In the case where $K\in \mathcal{P}_A$ is selfdual in the Tannakian sense (cf. Prop.~\ref{properties2}), 
the group $G_K$ can be interpreted as the largest subgroup of $\GL[H^0(A,K)^{\otimes n}]$  which stabilizes all subspaces $V\leq H^0(A,K)^{\otimes n}$ of the form $V=H^0(A,L),$ for $L$ a subquotient of the $n$-fold self-convolution $K^{\ast n}=K\ast \cdots \ast K$, hereby fixing elementwise subspaces of the form $V=H^0(A,\delta_0)$. When $A=E$ is an elliptic curve, the restriction of $K^{\ast n}$ to a suitable dense open subset of its support is a shift of a local system whose monodromy governs its decomposition into subquotients to a large part: \emph{the knowledge of the monodromy of $K\ast K,\, K\ast K\ast K,\ldots$ is therefore a crucial ingredient in the determination of the Tannakian group $G_K$}.


\bigskip

Let $E$ be an elliptic curve over an algebraically close field $\bar k$ and $\ell\neq \char(\bar k)$. The aim of this paper is to identify the Tannaka group $G_K$ for $K\in\PPE$ belonging to a certain class of self-dual unipotent \emph{seven-point $\QQ_{\ell}$-sheaves} (see Def.~\ref{def:7shef}), a generalization of perverse sheaves on $E$ obtained via the Beauville classification of families of elliptic curves \cite{Beauville}. We establish that \emph{for $K$ a seven-point sheaf on $E$, $G_K$ is isomorphic to the exceptional group $G_2$ (Thm.~\ref{thmapp})}.

This extends a geometric result first obtained by Katz in \cite{KatzEll} in the case of the Beauville classification: while Katz' proof relies on the evaluation of Frobenius traces for sufficiently general primes, our approach relies on the tensor decomposition of the monodromy of the $2$- and $3$-fold self-convolution of $K$ and on their rank properties, which we show is sufficient to identify $G_2$ (see Prop.~\ref{g2criterion} and \ref{propm3}). The summand-decompositions of $K\ast K$ and $K\ast K\ast K$ bring moreover some insight to two questions of Katz, see Rem.~\ref{rem:Katz}.

The difficulty in realizing $G_2$ as motivic Galois group -- by a classical Mumford-Tate argument, it can not be associated to the Absolute Hodge Cycle Motive of an abelian variety -- is also what makes this result of a particular interest (see also \cite{DR10}).

\bigskip

The plan of this paper is as follows: After some recollection on the Tannakian category $(\PPE,\ast)$ and after establishing some Tannakian properties of the convolution product in Section \ref{secrev}, we compute for $K_1,\ K_2\in \Perv(E)$ irreducible and not translation invariant perverse sheaves on $E$ the monodromy of the restriction of $K_1\ast K_2$ to its smooth locus in terms of elliptic braid groups in Section \ref{secbasicfibrations}. We then extend to this context a previous work of one of the authors with Wewers \cite{dw} on the variation of parabolic cohomology groups which gives an algorithm for this aim in Thm.~\ref{etalem}. This result, with further arguments, is applied to establish our main result Thm.~\ref{thmapp} in Section \ref{secbeauville}. The final argument relies on Magma computations whose code is given in Appendices \ref{app:braids} and \ref{app:MonMagma}.

\bigskip

The extension of the main result to the mixed case via variation of motives with supports -- thus leading to \emph{families of mixed motives on elliptic curves} --, and as announced in \cite{OWRD18} its application to Regular Inverse Galois Theory \emph{beyond the rigidity barrier} (see ibid.) will be the goals of forthcoming works.

\bigskip

\noindent\emph{Convention.} Throughout the article, multiplication of paths will  start from the left
(meaning that if we have a path-product $\alpha\beta$ then one first moves along $\alpha $ and then along $\beta$).

\bigskip

\noindent\emph{Acknowledgements.} While writing this paper, W.~S. was supported by Dr. Max R\"{o}ssler, the Walter Haefner Foundation, the ETH Z\"urich Foundation, and the Clay Foundation; B.C. was supported by DFG programme DE 1442/5-1 Bayreuth and DFG grant SPP1786.

\newpage
\section{Sheaf Convolution on Elliptic Curves}\label{secrev} 
Let $E$ be an elliptic curve over an algebraically closed field $\bar{k}$ and denote by $ a\colon E\times E\to E$ the addition map. Let $D^b_c(E,\bQl)$ denote the triangulated category of complexes of \'etales $\bQl$-sheaves over $E$. It contains the abelian subcategory $\Perv(E)$ of perverse sheaves -- complexes such that themselves and their Verdier dual are semi-perverse, i.e. satisfy a support-dimension condition, see \cite{BBD}, \S 4 -- as well as the category of constructible sheaves ${\rm Constr}(E),$ the latter placed in cohomological degree zero. Every $K\in D^b_c(E,\bQl)$ comes with a Verdier dual $D(K)\in D^b_c(E,\bQl)$ -- we refer to op. cit. for these generalities.

\subsection{Additive Convolution}
The {\it  convolution}
of two objects  $K_1,K_2 \in D^b_c(E,\bQl)$ is defined as 
\[ 
K_1\ast  K_2 = R^1a_*(\pr^*_1 K_1\otimes \pr^*_2K_2),
\] 
where $\pr_i:E\times E\to E$ denotes the $i$-th projection ($i=1,2$). 

\medskip

Consider the coordinate change $E\times E\to E\times E, (x,t)\mapsto (x,y:=x+t)$ 
so that with 
$d:E\times E\to E,(x,y)\mapsto y-x$ one has    
\begin{equation}\label{eq0}   
K_1\ast  K_2 = R^1\pr_{2*}(\pr^*_1 K_1\otimes d^*K_2).
\end{equation}
One therefore has the following expression for the stalk of $  K_1\ast  K_2$ at $y_0$:
\[ 
(K_1\ast  K_2)_{y_0}\simeq H^*(E,K_1\otimes K_2(y_0-x)),
\]	
where $K_2(y_0-x)=[x\mapsto y_0-x]^*K_2.$

\medskip

Let us collect some well-known properties of the convolution over an algebraically closed field (cf. \cite{WEI06}, \S 2.1 and \cite{KraemerWeissauer}, Prop. 10.1(b) and \cite{KRAPHD}, \S 2.1 for ii)).

\begin{prop}\label{properties}
	The triangulated category $(D^b_c(E,\bQl),\ast)$ is symmetric monoidal -- e.g. $K_1\ast K_2\simeq K_2\ast K_1$ and $(K_1\ast K_2)\ast K_3\simeq K_1\ast (K_2\ast K_3)$ -- more precisely:
	\begin{enumerate}[i)]
		\item The $\ast$-unit object is the skyscraper sheaf $\delta_0$ of rank one with support at the origin;
		\item Any $\ast$-invertible object is a skyscraper sheaf $\delta_x$ of rank one supported at some point $x\in E$;
	\end{enumerate}	
	Moreover, for $K_1,\ K_2\in D^b_c(E,\bQl)$:
	\begin{enumerate}[resume, label=\roman*)]
		\item Convolution commutes with shifts of complexes: 
		\[
		K_1[m]\ast K_2[n]=(K_1\ast K_2)[m+n];
		\]
		\item Convolution commutes with Verdier duality:
		\[
		D(K_1\ast K_2) \simeq D(K_1) \ast D(K_2).
		\]
	\end{enumerate}
\end{prop}

\subsection{A Perverse Tannakian Category}

As established for the middle product convolution over $\PP^1\setminus\{0,1,\infty\}$ by Deligne in \cite{DELRAC} Prop.~2, the role of the following property $P$ is to ensure for a certain category $\PPE$ of perverse sheaves to be stable under the convolution product, which then turns into a neutral Tannakian category. The situation is similar for abelian varieties, see \cite{KRA14} Thm.~3.8; we review this property for elliptic curves. 

\begin{defn} An object $K\in D^c_b(E,\bQl)$ has property $P$ if for all smooth rank-one sheaves 
	$L$ on $E$ the following holds:
	$$ H^i(E,K\otimes L)=0\quad \textrm{for}\quad i\neq 0.$$
\end{defn}

In this case, the following holds (cf.~\cite{KatzEll}, Lem.~2.1 and Cor~2.2, and \cite{KraemerWeissauer}, Thm.~7.1):

\begin{prop}\label{properties2} If $K \in D^c_b(E,\bQl)$ has property $P$, then $K\in \Perv(E)$. Moreover, the category $(\PPE,\ast)$ of semisimple perverse sheaves with property $P$ forms a Tannakian category neutralized by the fibre functor $\omega\colon K\mapsto H^0(E,K)$, with the skyscraper sheaf $\delta_0$ as identity, with dual 
	\[ 
	K \mapsto K^\vee:= [P\mapsto -P]^*DK
	\]
	and Tannakian dimension 
	\[
	\dim(K):=h^0(E,K)=\chi(K),
	\]
	where $\chi(K)$ denotes the Euler characteristic of $K.$ 
\end{prop}

For a perverse sheaf $K\in \PPE$, one forms the subtannakian category $\mathcal{P}_{E}(K)=\langle K\rangle$ of $\PPE$ that is finitely tensored generated by $K$, that is composed of subquotients of convolutions powers $(K\oplus K^{\vee})^{\ast n}$. The general Tannakian formalism gives an algebraic group $G_K$ called \emph{the Tannakian group of $K$} -- we refer to \cite{KRA14} \S 5 and Cor.~5.3 for details.			

\begin{rem}	\mbox{}	\begin{enumerate}
		\item Working in the category $\PPE$ should be seen as working in the \emph{quotient} category $\Perv(E)/Neg,$ where $Neg$ denotes the Serre subcategory of negligible objects of $\Perv(E)$ formed by the smooth objects on $E$ -- see \cite{KatzEll} Rem~2.3 and \cite{KRA14} \S 4 for general formalism and results.
		\item The Tannakian group $G_K$ attached to an irreducible object $K$ of $\PPE$ has the following useful characterization:
		it is the stabilizer inside $\GL(H^0(E,K))$ of all decompositions of 
		$H^0(E,K^{*n}*(K^\vee)^m)=H^0(E,K)^{\otimes n}\otimes H^0(E,K^\vee)^{\otimes m}\, (n,m \in \NN)$ which are induced 
		by the decomposition into irreducible subfactors of $K^{*n}*(K^\vee)^{*m}$, hereby fixing elementwise the contributions coming from $H^0(E,L)$ with $L=\delta_0$ a subfactor of $K^{*n}*(K^\vee)^{*m}$, cf. \cite{DELMIL82} Chapter II.
	\end{enumerate}
\end{rem}

We give three lemmata for  perverse sheaves in $\PPE$ that reflect their Tannakian properties. They were obtained by Katz as unpublished side results to \cite{KatzEll}, and included in lectures attended by one of the authors (W.~S.) in Princeton.

\medskip

\begin{lem}\label{WS1} Let $K_1$ and $K_2$ be perverse sheaves in $\PPE$. Then the generic rank 
	of $K_1 \ast  K_2$ is given as \[ \rk( K_1 \ast  K_2 ) = \rk(K_1)\chi(K_2) + \chi(K_1) \rk(K_2) .\] \end{lem} 

While a similar and more general statement exists for arbitrary abelian varieties in characteristic zero \cite{KRA16}, we provide a proof in our context and for the reader's convenience.

\begin{proof}  Let $\mathbf{x}_1=\{x_1,\dots,x_n\}$ be the singularities of $K_1$ and $\mathbf{x}_2=\{y_1,\dots,y_m\}$ be the singularities of $K_2$. Because $K_1\ast  K_2$ is a perverse sheaf, its generic rank is equal to minus its Euler characteristic at a generic point. We will calculate its Euler characteristic at any point $y \not\in \{ x_i+y_j | 1\leq i \leq n, 1 \leq j \leq m \}$. 
	At such a point, its Euler characteristic is $\chi (E , K_1 \otimes K_2(y-x))$.  The Euler characteristic of a complex of sheaves on an elliptic curves is a sum of local contributions at the singular points. (There would also be a global contribution, but the Euler characteristic of an elliptic curve is zero, so it vanishes). Because $K_2(y-x)$ is smooth in degree $-1$ of rank $\rk(K_2)$ in a neighbourhood of $x_i\in \mathbf{x}_1$, the local contribution to the Euler characteristic of $K_1 \otimes K_2(y-x)$ at the point $x_i$ is $-\rk(K_2)$ times the local contribution to the Euler characteristic of $K$ at $x_i$.  So the total contribution at $x_1,\dots,x_n$ is $-\chi(K_1)\rk(K_2)$. Similarly, $K_1$ is smooth in degree $-1$ of rank $\rk(K_1)$ in a neighbourhood of $y-y_j$ for $y_j\in\mathbf{x}_2$, so the local contribution to the Euler characteristic of $K_1 \otimes K_2(y-x) $ at the point $y-y_j$ is $-\rk(K_2)$ times the local contribution to the Euler characteristic of $K_2(y-x) $ at $y-y_j$, which is the local contribution to the Euler characteristic of $K_2$ at $y_j$.  Hence the total Euler characteristic of $K_1 \otimes K_2(y-x) $ is $-\rk(K_1)\chi(K_2) -\chi(K_1) \rk(K_2) $, as desired.
\end{proof}

\begin{lem}\label{WS2} Let $K\in \PPE$ be a perverse sheaf whose associated representation of the Tannakian group $G_K$ has finite order. Then $K$ is a sum of skyscraper sheaves on points of $E$ of finite order. 
\end{lem}

We provide a proof for this result which is different from Katz' one.
(For a proof for an arbitrary abelian variety see \cite{Wei12}, Thm.~3, and for Katz's argument in a different setting see \cite{KatzMellin}, Thm.~6.4.) 

\begin{proof}  The $n$th convolution power $K^{\ast n}$ is a sum of irreducible perverse sheaves corresponding 
	to the irreducible representations of the finite Tannakian group of $K$.
	Each appears with multiplicity at most the Tannakian dimension of $K^{\ast n}$, which is $\chi(K)^n$, and has some finite rank,
	so the total rank of $K^{\ast n}$ is $O(\chi(K)^n)$. By inductively applying Lem.~\ref{WS1}, 
	we see that  the rank of $K^{\ast n} $ is $n \chi(K)^{n-1} \rk(K)$. It follows that $\rk(K)=0$, so $K$ is a sum of skyscraper sheaves.
	Because skyscraper sheaves at points of infinite order have monodromy group $\mathbb G_m$, 
	it follows that all these are at points of finite order.
\end{proof}

\begin{lem}\label{WS3} Let $K\in \PPE$ be an irreducible perverse sheaf. If the associated representation of the Tannakian group $G_K$ is not irreducible when restricted to its identity component, then $K$ is isomorphic to its own translate by some nontrivial point of $E$. \end{lem}

\begin{proof} Let $V$ be the representation associated to $K$, and $G_0$ the identity component of its Tannakian group $G_K$. Because $V$ is irreducible as a representation of $G_K$, it is semisimple as a representation of $G_0$, hence if $V$ is not irreducible as a representation of $G_0$ then $\Hom_{G_0}(V,V)$ has dimension greater than one. The $G_K$-invariant subspace of $\Hom_{G_0}(V,V)$ is equal to $\Hom_{G_K}(V,V)$, which has dimension one, so $\Hom_{G_0}(V,V)$ is a nontrivial representation of $G_K$. 
	
	Because $\Hom_{G_0}(V,V)$ factors through $G_K/G_0$, it has finite monodromy, and thus by Lem.~\ref{WS2} it is a sum of one-dimensional objects corresponding to skyscraper sheaves. For each of these one-dimensional representations $W$, there is a nontrivial map $V \otimes W \to V$, hence an isomorphism as $W$ is one-dimensional, so an isomorphism $K \ast  \delta_a \cong K$, which implies that $K$ is isomorphic to its own translate by $a\in E$. 
\end{proof}

Notice that the same definition of convolution product holds for topological constructible sheaves with respect to the complex topology. Via analytification and the \'etale-singular cohomology comparison isomorphism, this reduces the study of the geometric properties of the convolution to the topological context (cf.~\cite{dw}, Thm.~3.2).

Let thus $K=j_*L[1]\in \PPE,$ where $L$ is a local system on a dense open subset $j:U\hookrightarrow E$, which is thus considered over $\CC$ and relatively to the complex analytic topology. Let $\x=\{x_1,\ldots,x_p\}:=E\setminus U.$ 
Recall that the fundamental group 
of $E\setminus \x$ has the following presentation:
\[
\pi_1(E\setminus \x,x_0)=\langle \gamma_1,\ldots, \gamma_r, \alpha, \beta\mid \gamma_1\cdots \gamma_r[\alpha,\beta]=1\rangle,
\]
where $\gamma_i\, (i=1,\ldots,r)$ is the homotopy class of a  simple closed loop going counterclockwise 
around the missing point $x_i,$ where 
$\alpha,\beta$ are the homotopy classes the usual generators of $\pi_1(E)$ and where
$[\alpha,\beta]=\alpha\beta\alpha^{-1}\beta^{-1}$ (cf.~Fig.~\ref{fig:Braid_Conf2}).
Therefore, the sheaf $L$ (and hence $K=j_*L[1]$, by the uniqueness of the intermediate extension) 
corresponds, via its monodromy representation, to 
its {\it monodromy tuple}
$$ T_L=T_K:=(A_1,\ldots,A_r,A,B)\in \GL_n(\bQl)$$ 
satisfying 
$$ A_1\cdots A_r[A,B]=1,$$ where 
$A_i\,(i=1,\ldots,r)$ is the image of $\gamma_i$ in $\GL(L_{x_0})\simeq \GL_n(\bQl)$  and where 
$A,B$ are the images of $\alpha$ and $\beta$ (resp.). This gives the following explicit formula for the Tannakian dimension for irreducible intermediate extensions in $\PPE:$		

\begin{prop}
	Let $L\in \PPE$ be an irreducible nonconstant smooth sheaf of rank $n$ on $E\setminus \x$ and let $K:=j_*L[1]$, then $h^0(E, K)=h^1(E,j_*L)$ and 
	\[
	\dim(K)=\chi(K)= rn-\sum_{i=1}^r\dim(V^{ A_i }),
	\]
	where $V=\bQl^n$ and where $V^{A_i}$ denotes the fixed space of $A_i.$ 
\end{prop}	

\begin{proof} The result is straightforward by the additivity of the Euler characteristic
	\[
	\chi(E,j_*L)=\chi(E\setminus \x,L)+\chi(\x,L|_{\x})=rn-\sum_{i=1}^r h^0(x_i,L|_{x_i}).
	\]
\end{proof}

We work from now on over $\CC$, in the analytic context, and we view the objects relatively to the complex analytic topology. 


\section{Monodromy of Convolution and Elliptic Braid Groups}\label{secbasicfibrations}

Let $K_1,\ K_2$ be two intermediate extensions of smooth irreducible perverse sheaves $L_1[1],\ L_2[1]$ contained in $\PPE$. Denote by $U_i\subset E$, $i=1,2$, their smooth locus, and by respectively $\x_1=\{x_1,\ldots,x_p\}$ and $\x_2=\{y_1,\ldots,y_q\}$ the singular set of $K_1$ and $K_2$. The smooth locus of their convolution $K_1\ast K_2$ defined in Eq.~\eqref{eq0} is given by:
\[
\mathbb{V}:=\pr^{-1}_1(U_1)\cap d^{-1}(U_2)\cap \pr^{-1}_2(E_y\setminus \x_1\ast \x_2) \subseteq E_x\times E_y
\]
where $\x_1\ast \x_2:=\{x_i+y_j\mid i=1,\ldots,p,\, j=1,\ldots,q\}$. 

\medskip

Since the restriction of $\pr_2$ to $\mathbb{V}$ induces a locally trivial fibration $\pr_2 \colon \mathbb{V} \to E_y\setminus \x_1\ast \x_2$ -- whose fibres are copies of $E$ with $p+q$ points deleted -- one obtains:
\begin{equation}\label{eq:EllFib}
1\to \pi_1(E\setminus (\x_1\cup y_0-\x_2),x_0) \to \pi_1(\mathbb{V},(x_0,y_0)) \to \pi_1(E_y\setminus \x_1\ast \x_2,y_0) \to 1. 
\end{equation}
A splitting of the sequence above gives rise to an action of the fundamental group of the base on the fiber, which via the equivalence of categories between local systems over $E$ and $\pi_1(E,y_0)$-representations encodes the variation of the local system associated to $K_1\ast K_2$ over $E\setminus \x_1\ast\x_2$, see \S \ref{sec:MonConv} for details.

	\medskip

We compute this monodromy in terms of a Birman split-fibration that identifies $\pi_1(E\setminus (\x_1\cup y_0-\x_2),x_0) $ and $\pi_1(E_y\setminus \x_1\ast \x_2,y_0)$ with elliptic braids groups.

\subsection{Elliptic Braids, Presentation and Fibration}\label{sub:EllBr}
For $n\in \NN_{>0},$  let $F_n(E)=E^n\setminus \Delta$ denote the configuration space of $n$ points on an elliptic surface whose homotopy group is \emph{the pure elliptic surface braid group $P(n,E):=\pi_1(F_n(E),(x_1,\ldots,x_n))$}, cf.~\cite{Bellingeri}, 
\S~2.1 and \cite{GuaschiGoncalves}. In what follows, we consider the presentation of $P(n,E)$ as given in \cite{Bellingeri}, Thm.~5.1 within the full elliptic braid group $\widetilde{P}(n,E)$ by generators:
\[
\widetilde{P}(n,E)=\langle \sigma_1,\ldots,\sigma_{n-1},\alpha,\beta\rangle
\]				
and the relations of the following kind:

\begin{itemize}
	\item Braid relations:
	\begin{equation}\tag{B1-B2}
	\sigma_i\sigma_{i+1}\sigma_i=\sigma_{i+1}\sigma_i\sigma_{i+1}\quad \textrm{ and } \quad \sigma_{i}\sigma_j=\sigma_j\sigma_i \textrm{ for } |i-j|\geq 2; 				
	\end{equation}
	
	\item Mixed relations:
	\begin{gather}
	[\alpha,\sigma_i]=[\beta,\sigma_i]=1\quad \textrm{ for } i\neq n-1\tag{M1}\\ 
	\sigma_{n-1}^{-1}\alpha\sigma_{n-1}^{-1}\alpha=\alpha\sigma_{n-1}^{-1}\alpha\sigma_{n-1}^{-1}\quad \text{ and }\quad  \sigma_{n-1}^{-1}\beta\sigma_{n-1}^{-1}\beta=\beta\sigma_{n-1}^{-1}\beta\sigma_{n-1}^{-1} \tag{M2-M3}\\
	\sigma_{n-1}^{-1}\alpha\sigma_{n-1}^{-1}\beta=\beta\sigma_{n-1}^{-1}\alpha\sigma_{n-1} \tag{M4}\\
	[\alpha,\beta]=\sigma_{n-1}\sigma_{n-2}\cdots\sigma_{1}^2\cdots \sigma_{n-2}\sigma_{n-1}\tag{M5}
	\end{gather}
\end{itemize}
where $[g_1,g_2]:=g_1g_2g_1^{-1}g_2^{-1}$. Following loc. cit. $P(n,E)$ is more precisely generated as a subgroup of $\widetilde{P}(n,E)$ by:
\begin{itemize}
	\item \emph{Global braids:} 
	\begin{equation}\tag{P1}
	\alpha^{\sigma_{n-1}\cdots \sigma_k} \textrm{ and } \beta^{\sigma_{n-1}\cdots \sigma_k} \text{ for } k=1,\ldots,n-1
	\end{equation}
	\item \emph{Local braids:} 
	\begin{equation}\tag{P2}\label{eq:P2}
	\beta_{i,j}:=\beta_{j-1}\cdots \beta_{i+1}\beta_{i}^2\beta_{i+1}^{-1}\cdots \beta_{j-1}^{-1} \text{ for } 1\leq i<j\leq n
	\end{equation}
\end{itemize}
where we write $\beta_i:=\sigma_i^{-1}$.

\medskip

The forgetting map
\begin{equation}\label{eq2}  
\theta:F^{n+1}(E)\to F^{n}(E), (x_1,\ldots,x_n,x_0)\mapsto (x_1,\ldots, x_{n})
\end{equation} 
is a locally trivial fibration whose fibres are again copies of $E$ with $n$ points deleted. The long exact homotopy sequence gives rise to a short exact sequence:
\[
1\to \pi_1(E\setminus \{x_1\ldots,x_{n}\},x_0)\to P(n+1,E)\to P(n,E)\to 1,
\] 
which by \cite{BIR69}, Proof of Thm.~5 admits an algebraic splitting.

\subsection{Elliptic Braids and Basic Fibration} For $n=p+q$ the homotopy elliptic fibration \eqref{eq:EllFib} and the homotopy sequence above are then put together via the morphism of fibrations:
\begin{equation}\label{splitdiag0}\begin{CD} \mathbb{V}
@>>>  E_y\setminus \x_1\ast \x_2 \\
@V{\tilde{\phi}}VV   @V{\phi}VV    \\
F^{p+q+1}(E) @>>> F^{p+q}(E)\,  ,  \\
\end{CD}\end{equation}
where $$\tilde{\phi}(x,y)=(x_1,\ldots,x_p,y-y_1,\ldots y-y_q,x)\textrm{ and } \phi(y)=(x_1,\ldots,x_p,y-y_1,\ldots y-y_q).$$

The homotopy sequences of the rows give rise to the following commutative diagram of homotopy sequences:
\begin{equation}\label{splitdiag}
\begin{gathered}
\xymatrix@C=1.4em{
	1\ar[r] & \pi_1(E\setminus (\x_1\cup y_0-\x_2),x_0)\ar[r]\ar[d]^= & \pi_1(\mathbb{V},(x_0,y_0))\ar[r]\ar[d]^{\tilde{\phi}_*} & \pi_1(E_y\setminus \x_1\ast \x_2,y_0)\ar[r]\ar[d]^{\phi_*}& 1\\
	1\ar[r] & \pi_1(E\setminus (\x_1\cup y_0-\x_2),x_0)\ar[r] & P(p+q+1,E)\ar[r] & P(p+q,E)\ar[r] & 1\, .
}
\end{gathered}
\end{equation}
The Birman algebraic splitting of the lower row induces an action of $\pi_1(E_y\setminus \x_1\ast \x_2,y_0)$ 
on $\pi_1(E\setminus (\x_1\cup y_0-\x_2),x_0)$ via its image under the fundamental morphism $\phi_*:$
\begin{equation}\label{eq:FundPhi}
\phi_*\colon\pi_1(E\setminus \x_1\ast \x_2,y_0) \to P(p+q,E)
\end{equation}

\medskip

\begin{figure}[h!]
	\centering
	\begin{subfigure}[b]{0.40\textwidth}
		\includegraphics[width=0.95\linewidth]{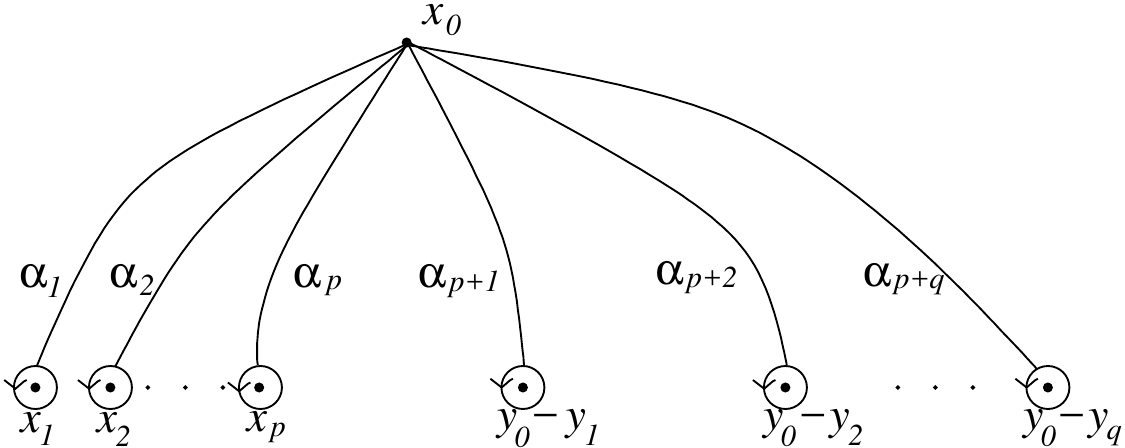}
		\caption{Braids $\alpha_{i}$.}\label{fig:Pure_Gen}
	\end{subfigure}
	\begin{subfigure}[b]{0.54\textwidth}
		\includegraphics[width=0.95\linewidth]{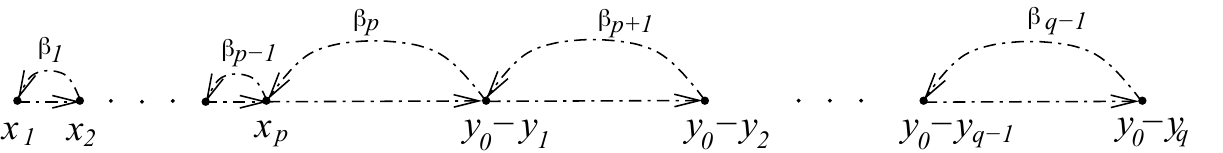}
		\caption{Braids $\beta_{i}$.}\label{fig:Loc_Braids}
	\end{subfigure}
	\caption{Braids on $E\setminus (\x\cup y_0-\y)$}
	\label{fig:Braid_Torus}
\end{figure}

According to the presentation of $P(E,n)$ of \S \ref{sub:EllBr} and the commutative diagram~\eqref{splitdiag}, one obtains more precisely:
\begin{prop}\label{prop:presPfibre}
	The fundamental group 
	$ \pi_1(E\setminus (\x_1\cup y_0-\x_2),x_0)$ admits an action of the base $\pi_1(E_y\setminus \x_1\ast \x_2,y_0)$, and has a presentation
	$$ \pi_1(E\setminus (\x_1\cup y_0-\x_2),x_0)=\langle \alpha_1,\ldots,\alpha_n,\alpha,\beta \mid 
	\alpha_1\cdots \alpha_n[\alpha, \beta]=1 \rangle$$
	where $ \alpha_1:=\beta_{1,n+1},\ldots,\alpha_{n}:=\beta_{n,n+1}$, and where the $\alpha_i$ and the $\beta_{i,j}$ are given in \eqref{eq:P2} -- cf. Fig.~\ref{fig:Braid_Torus}.
\end{prop}	

\medskip

The following section computes the morphism $\phi_*$ of Eq. \eqref{eq:FundPhi} and the previous action on the fiber in terms of $P(n,E)$.

\subsection{Computing the Monodromy of Convolution}\label{secellbraids}
We now compute the action  of $\pi_1(E_y\setminus \x_1\ast \x_2,y_0)$ on $\pi_1(E\setminus (\x_1\cup y_0-\x_2),x_0)$ via the fundamental morphism $\phi_*$ of Eq. \ref{eq:FundPhi} and using the presentation of Prop.~\ref{prop:presPfibre} in $P(n,E)$. We first compute the image by $\phi_*$ of the \emph{local braids} $\delta_{i,j}$ around the $pq$ ramification points $\x_1\ast \x_2$, and their action, then whose of the two \emph{global braids} $\hat{\alpha}$ and $\hat{\beta}$ around the handle. 

\medskip

\begin{figure}[h!]
	\centering
	\includegraphics[width=0.7\textwidth]{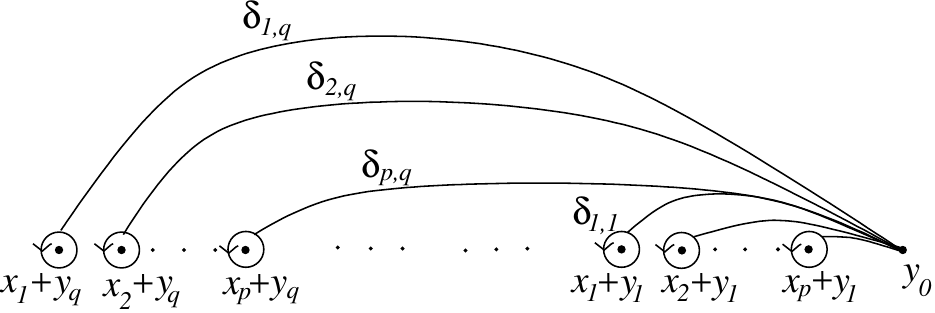}
	\caption{Local Braids on $E\setminus \x_1* \x_2$.}
	\label{fig:LocBrBase}
\end{figure}

\medskip

We choose the \emph{local braids} of the base $E\setminus \x_1* \x_2$ as in Fig.~\ref{fig:LocBrBase}. Then one has the following result:
\begin{prop}\label{propdeltt} Assume that the singular loci $\x_1\ast\x_2$ is generic, i.e.
	$\#(\x_1\ast \x_2)=pq.$ The morphism $\phi:\pi_1(E\setminus \x_1\ast \x_2,y_0) \to P(p+q,E) $ of Eq.~\ref{eq:FundPhi} is given on the local braids by:
	\begin{subequations}
		\begin{align}
		\phi(\delta_{i,1}) &=\beta_{i,p+1},\; i=1,\ldots,p,\label{eqimageofdeltai} \\
		\phi(\delta_{i,j}) &=\beta_{i,p+1}^{\beta_{p+1}\cdots \beta_{p+j-1}},\; i=1,\ldots,p,\; j=2,\ldots,q \,.\label{eqimageofdeltai2} 
		\end{align}
	\end{subequations}
\end{prop}

\begin{proof} First recall that $E$ is represented by $\RR^2/\Lambda,$  with $\Lambda$ a lattice in $\RR^2.$ Hence, locally, the addition in $E$ 
	is given by the addition in $\RR^2.$ Therefore, if the singular loci
	$\x_1, \x_2, y_0-\x_2$ are assumed to be closely together in the analytic topology then, locally, one recovers the arrangement of singular points of 
	$L_1\boxtimes L_2$ of the additive 
	convolution as in Fig.~\ref{fig:ConvG0}.
	
	\begin{figure}[h!]
		\centering
		\includegraphics[width=0.7\textwidth]{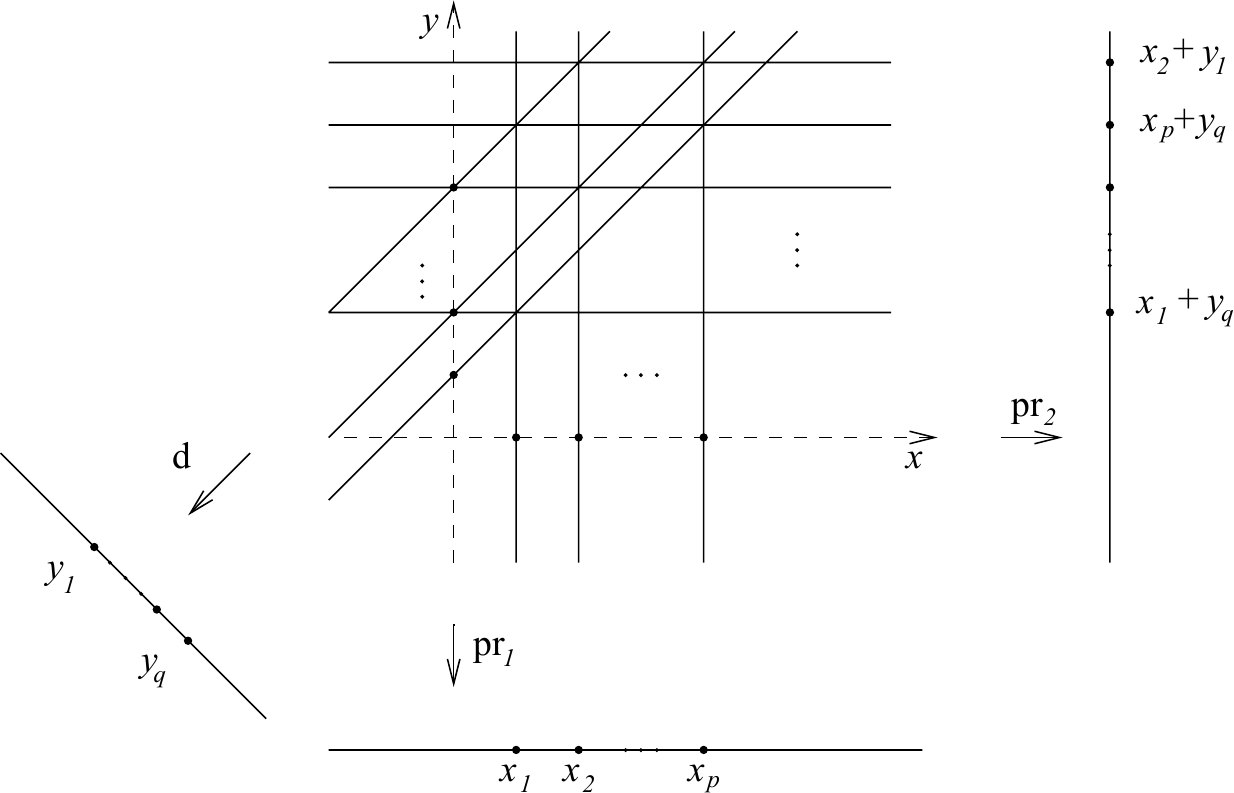}
		\caption{Local Convolution $E\setminus \x_1* \x_2$.}
		\label{fig:ConvG0}
	\end{figure}
	
	\medskip
	
	The formulae are then direct consequences of \cite{DJ18}, Prop.~3.2.3, using that the restriction of singularity arrangement considered there to a disk containing $\x_1\ast \x_2$ is 
	the same as the above arrangement.
\end{proof}

Recall that $\beta_i\, (i=1,\ldots,n-1)$ acts by componentwise conjugation  with $\beta_i^{-1}$ on $(\alpha_1,\ldots,\alpha_n,\alpha,\beta)$ via 
\begin{equation}\label{eqbetai} 
\begin{split}
\beta_i(\alpha_1,\ldots,\beta)&=(\beta_i\alpha_1\beta_{i}^{-1},\ldots,\beta_i\beta\beta_i^{-1})\\
&=(\alpha_1,\ldots,\alpha_{i-1},\alpha_{i+1}, \alpha_{i+1}^{-1}\alpha_i\alpha_{i+1},\alpha_{i+1},\ldots,\alpha_n,\alpha,\beta),
\end{split}
\end{equation}
inducing an action from the right on local systems which transforms a monodromy tuple $(A_1,\ldots,A_n,A,B)$
into
$$ (A_1,\ldots,A_n,A,B)^{\beta_i}=(A_1,\ldots,A_{i+1},A_{i+1}^{-1}A_iA_{i+1},\ldots,A_n,A,B),$$
cf.~\cite{dw}, \S~2.2.
Using this rule, one obtains the operation of the above braids $\phi(\delta_{i,j})$ on $ (\alpha_1,\ldots,\alpha_n,\alpha,\beta),$ crucial for the 
results below (cf.~Lem~\ref{etalem}).

\medskip

To obtain this operation also for the \emph{global paths of the base}, we proceed as follows:
Let $\hat{\alpha},\hat{\beta}$ be the global elements of $\pi_1(E_y\setminus \x_1\ast \x_2,y_0)$ such that 
\[ \delta_{1,q}\delta_{2,q}\cdots \delta_{p-1,1}\delta_{p,1}[\hat{\alpha},\hat{\beta}]=1.\] 
The operation of $\phi(\hat{\alpha})$ and $\phi(\hat{\beta})$ on 
$\tau_{L_1\otimes L_2(y_0-x)}:=(\alpha_1,\ldots\alpha_{p+q},\alpha,\beta)$, which is defined componentwise by conjugation with respect to the semidirect product structure as in diagram~\eqref{splitdiag}, e.g.:
$$ \tau_{L_1\otimes L_2(y_0-x)}^{\phi(\hat{\alpha})}:=({\phi(\hat{\alpha})}^{-1}\alpha_1{\phi(\hat{\alpha})},\ldots,{\phi(\hat{\alpha})}^{-1}\alpha_{p+q}{\phi(\hat{\alpha})},{\phi(\hat{\alpha})}^{-1}\alpha{\phi(\hat{\alpha})},{\phi(\hat{\alpha})}^{-1}\beta{\phi(\hat{\alpha})})$$
is given by the following proposition.

\begin{prop}\label{propdeltt2} For $\hat{\alpha}$ and $\hat{\beta}$ as above, the operation of $\phi(\hat{\alpha})$ and $\phi(\hat{\beta})$ is given by:
	\[
	\begin{split}
	\tau_{L_1\otimes L_2(y_0-x)}^{\phi(\hat{\alpha})}=&
	(\alpha_1^{\alpha_{p+1}\cdots\alpha_{p+q}},\ldots, \alpha_p^{\alpha_{p+1}\cdots\alpha_{p+q}},\alpha_{p+1}^{\alpha^{-1}},\ldots,\alpha_{p+q}^{\alpha^{-1}},
	\alpha, \alpha_{p+q}^{-1}\cdots \alpha_{p+1}^{-1}\beta)\\
	\tau_{L_1\otimes L_2(y_0-x)}^{\phi(\hat{\beta})}=&
	(\alpha_1^{\beta\alpha_{p+q}^{-1}\cdots\alpha_{p+1}^{-1}\beta^{-1}},\ldots, \alpha_p^{\beta\alpha_{p+q}^{-1}\cdots\alpha_{p+1}^{-1}\beta^{-1}},
	\alpha_{p+1}^{\beta^{-1}},\ldots,\alpha_{p+q}^{\beta^{-1}},
	\alpha_{p+1}\cdots\alpha_{p+q}\alpha,\beta).
	\end{split}
	\]
\end{prop}

\begin{proof}	 \begin{figure}[h!]
		\centering \footnotesize
		\begin{subfigure}{0.32\textwidth}
			\def\svgwidth{0.95\linewidth}
			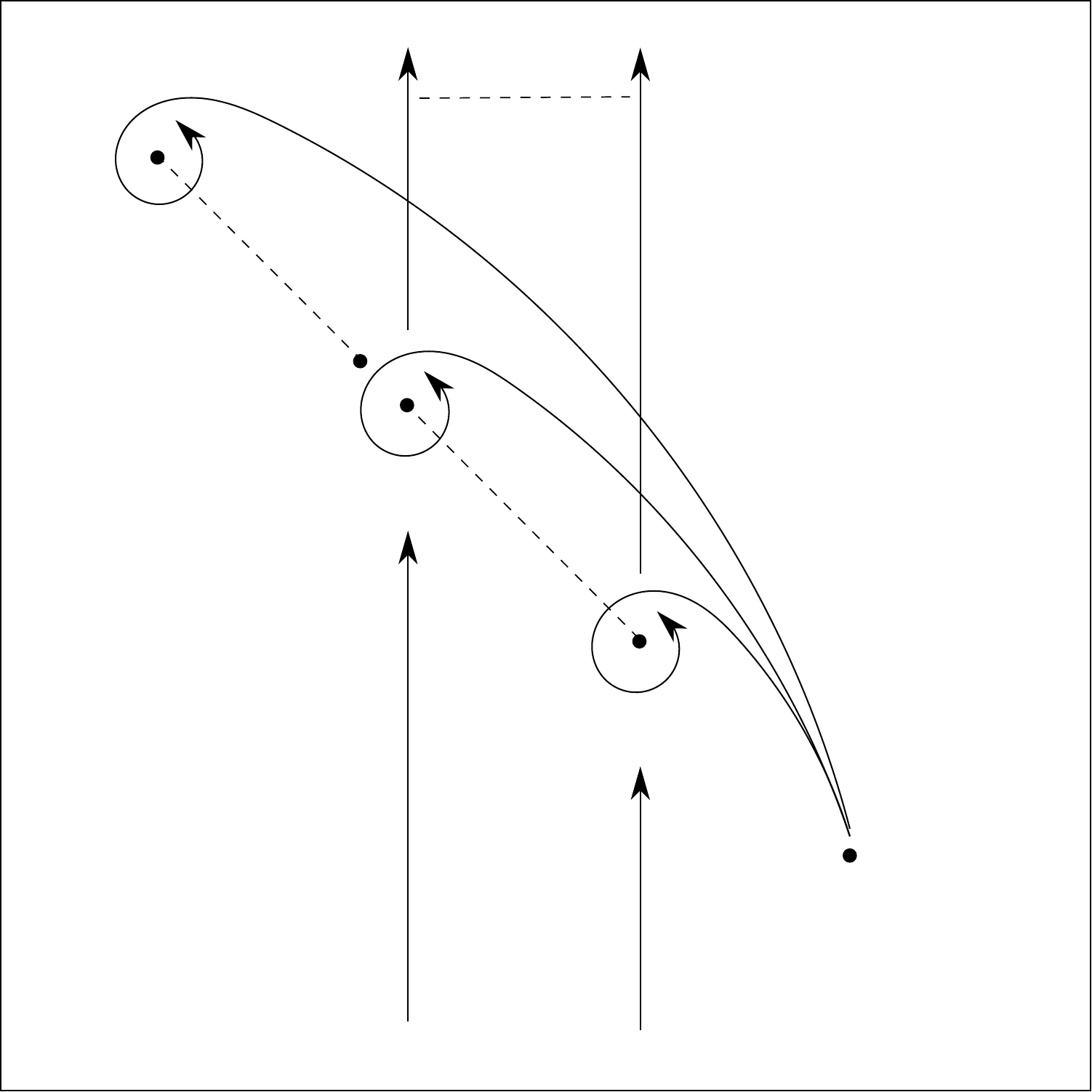{}
			\caption{Braids Configuration.}\label{fig:Braid_Conf}
		\end{subfigure}
		\begin{subfigure}{0.32\textwidth}
			\def\svgwidth{0.95\linewidth}
			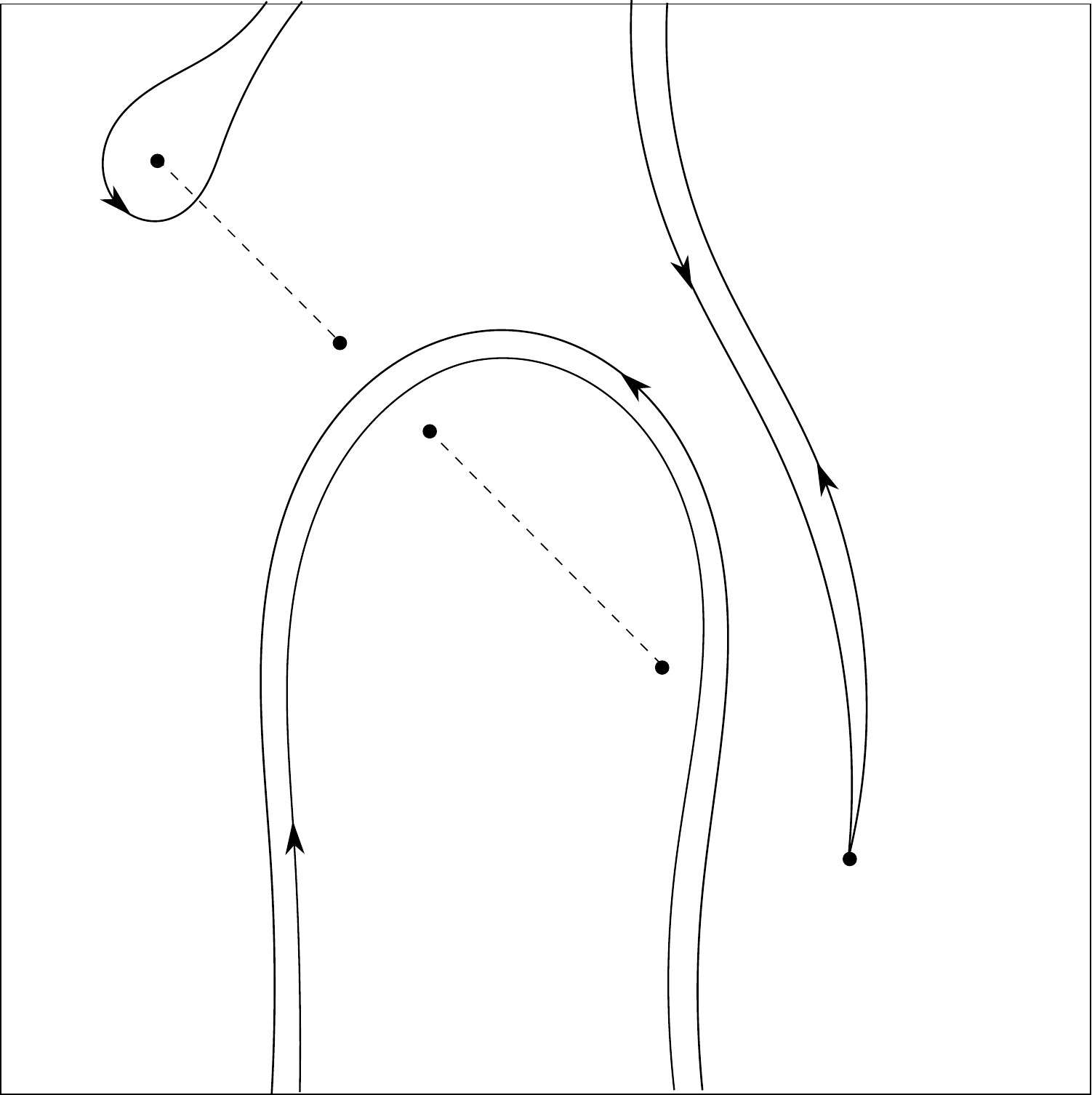{}
			\caption{Case $i=1$.}\label{fig:Braid_Act1}
		\end{subfigure}
		\begin{subfigure}{0.32\textwidth}
			\def\svgwidth{0.95\linewidth}
			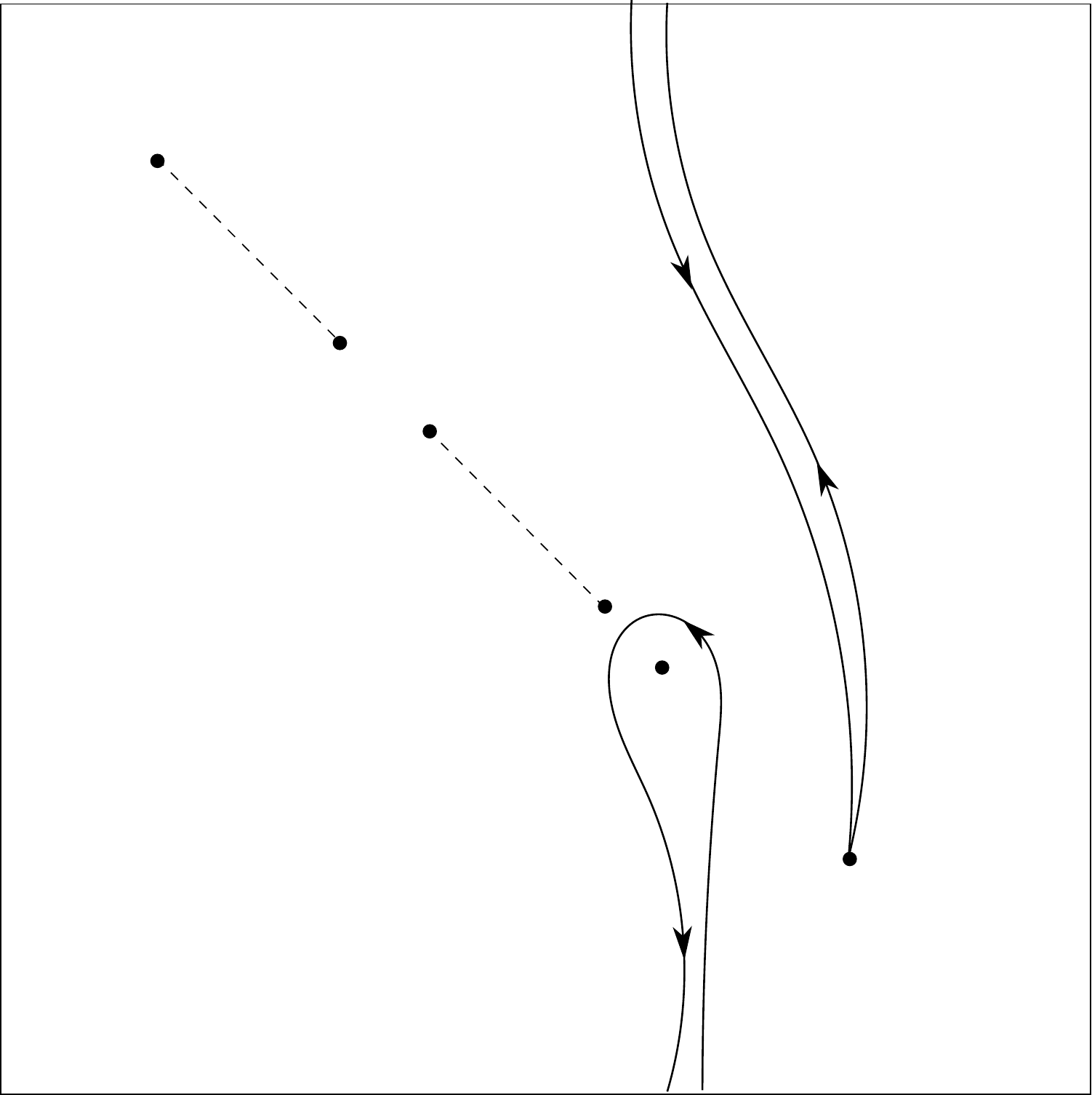{}
			\caption{Case $i=p+q$.}\label{fig:Braid_Act2}
		\end{subfigure}
		
		\caption{Action of $\phi(\hat{\beta})$ on $\alpha_i$}
		\label{fig:Braid_Acts}
	\end{figure} The computation involves two cases, depending whether $i\in \llbracket 1,p \rrbracket$ or $i\in \llbracket p+1,p+q \rrbracket$. For $\beta$ the proof is  straightforward by Fig.~\ref{fig:Braid_Acts} as in \cite{BIR69}, proof 3 of Corollary 5.1 (where the elliptic curve $E$ is represented by a square with opposite edges identified). 
	The action of $\alpha$ is obtained similarly.
\end{proof}

\section{The Monodromy Transformation of Cocycles of Convolutions}\label{sec:MonConv}
Let $\x_1\subseteq E$, $\x_2\subseteq E$ and let 
$L_1,\ L_2$ be two irreducible local systems respectively on $U_1:=E\setminus \x_1,\ U_2:=E\setminus \x_2$. We assume that the associated perverse intermediate extensions $K_1:=j_*L_1[1],\ K_2:=j_*L_2[1]$ are in $\Perv(E)$. By construction: 
\[ 
L_1\ast L_2:=\mathcal{H}^{-1}(K_1\ast K_2)|_{E\setminus \x_1\ast \x_2}
\] 
is a local system on $S=E\setminus \x_1\ast \x_2$ with monodromy representation:
\begin{equation}\label{eq:ConvMonRep} 
\eta\colon\pi_1(S)\to \GL(H^1(E,j_*(L_1\otimes L_2(y_0-x)))).
\end{equation}

\medskip

On the other hand, for $j\colon U\hookrightarrow E$ and $L$ sheaf over $U$ as before, consider the parabolic cohomology $H^1_p(U,L):=H^1(E,j_\star L)$. Since we deal with families of elliptic curves $E/S$ that are \emph{isotrivial} -- since given by moving singularities--, it follows from \cite{dw} \S 2.1, that the representation of Eq. \eqref{eq:ConvMonRep} can be computed in terms of the parabolic cohomology $H ^1_p(E\setminus (\x\cup y_0-\y),L_1\otimes L_2(y_0-x))$ -- see Eq. \eqref{splitdiag} for the fibration.

\medskip

After a short remark on Thom-Sebastiani isomorphism, we compute the monodromy representation $\eta$ via the identification of $H^1_p(U,L)$ to an explicit finitely presented algebra $W_T$ associated to the monodromy tuple $T$ of $L$ -- see \S \ref{sub:ParaCoc}

\subsection{Local Monodromy and Thom-Sebastiani}
An a priori knowledge of the local monodromy of $K_1\ast K_2$ is provided by Thom-Sebastiani Theorem that relates the vanishing cycle of $K_1\ast K_2$ to those of $K_1$ and $K_2$. 

\medskip

Recall the notion of vanishing cycles $\varphi(x_i)$ of an intermediate extension $K=j_*L[1]$, where $L$ is a local system on $E\setminus \x$ and $j\colon E\setminus \x\to E$ the open inclusion: this is a skyscraper sheaf supported on $x_i\in E$ with global sections isomorphic to ${\rm Im}(A_i-1)$ equipped with the induced operation of $A_i$ (resp. $\alpha_i$). Notice that the skyscraper sheaf structure is sometimes neglected and simply viewed as a vector space with an $A_i$-action. The local monodromy of $L$ at $x_i$ is uniquely determined by $\varphi_{x_i}(K)$.

Since the vanishing cycles of a skyscraper sheaf $\delta_{x_i}$ coincide by construction with the sheaf itself, this leads to the notion of \emph{vanishing cycles for semisimple perverse sheaves}. The local monodromy of $K_1\ast K_2$ is now determined by the analytical Thom-Sebastiani Theorem.

\begin{thm}[Thom-Sebastiani]\label{thmthomseb} Under the assumptions on $K_1,K_2$ made above, the convolution $K_1\ast K_2$ is semisimple	and one has a short exact sequence of perverse sheaves:
	\[
	0\to j_*j^*(K_1\ast K_2)\to K_1\ast K_2 \to \delta\to 0,			
	\]
	where $j:E\setminus \x_1\ast \x_2\hookrightarrow E$ and where $\delta$ is zero, unless 
	there exists $x_0\in E$ such that $K_2(x_0-x)\simeq D(K_1)$ and  $\delta\simeq \delta_{x_0}.$ In both cases one has an isomorphism:
	\begin{equation*}
	\varphi_{x_i+y_j}(K_1\ast K_2)\simeq \bigoplus_{(s,t)\,: \,x_s+y_t=x_i+y_j}\varphi_{x_s}(K_1)\otimes \varphi_{y_t}(K_2).
	\end{equation*}
\end{thm}

\begin{proof} See \cite{Weiss1} for the first claim. The second claim follows from the fact that, locally, the elliptic curve convolution looks like the additive convolution on $\AA^1$ (cf. the discussion in the proof of Prop.~\ref{propdeltt}) and Thom-Sebastiani Theorem over $\mathbb C$ as explained in \cite{ds}, proof of Thm.~3.2.3. 
\end{proof}

Notice that for $K\in \Perv(E)$ a $\bar{\QQ}_\ell$-perverse sheaf on $E_{\bar{k}}$, one forms the \emph{vanishing cycle} $\varphi_{x_i}(K)\in D^b_c(E_{x_i},\bQl)$ of $K$ at $x_i\in E_{\bar{k}}$ as in \cite[Exp. XIII]{SGA7}-- by an argument of Gabber, the vanishing cycle functor preserves the perversity, see \cite{BBD}, \S 4.4. A Thom-Sebastiani for the additive convolution product over $E$ can then be deduced from Thm.~4.5 and Rem.~11 (c) of \cite{ILL16}, and one recovers the result above by analytification.

\subsection{Parabolic Cohomology and Cocyles Algebra}\label{sub:ParaCoc}
As mentioned above, linear automorphisms act from the right, i.e., if $A
\in \GL(V)$ and $v\in V,$ then $vA$ denotes the image of $v$ under
$A.$ Let $L$ be a local system of rank $n$ on $U=E\setminus \x\ (\x=\{x_1,\ldots,x_r\})$ and $j_*L$ is extension to $E.$
As in \cite{dw}, Prop.~1.1(i), one concludes that $H^n(U,L)\simeq H^n(\pi_1(U),V),$
where
$V\simeq \bQl^n$ is the $\pi_1(U)$-module underlying the monodromy representation of $L.$
One further concludes as in loc. cit. that
\[
H^1(E,j_*L)=\im[H^1_c(U,L)\to H^1(U,L)].
\] 
As usual, the first cohomology group $H^1(\pi_1(U),V)$ is formed, 
up to equivalence induced by exact cycles, by maps
$\delta:\pi_1(U)\to V$ satisfying the {\it  cocycle rule }
\begin{equation}\label{cocyclerule} \delta(\gamma_1\gamma_2)=\delta(\gamma_1)\rho_L(\gamma_2)+\delta(\gamma_2)\quad (\gamma_1,\gamma_2\in \pi_1(U))\,\end{equation}
where $\rho_L:\pi_1(U)\to \GL(V)$ denotes the monodromy representation of $L.$ 
Moreover, under the isomorphism $H^1(U,L)\simeq H^1(\pi_1(U),V),$ a cocycle is in the image 
of $H^1_c(U,L)$ if and only if for each $i \in \{1,\ldots,r\}$ the image $v_i:=\delta(\alpha_i)$ lies in the image of 
$A_i-1,$ cf.~\cite{dw}, Lem.~1.2. The evaluation map 
$$ \delta \mapsto (\delta(\alpha_1),\ldots, \delta(\alpha_r),\delta(\alpha),\delta(\beta))\in V^{r+2}$$
and the above cocycle rule
gives the following explicit description of $ H^1(E,j_*L)$ in terms of the monodromy tuple
$T=T_L$  of $L:$ 
\begin{lem}\label{lemH1}
	There is  an isomorphism $ H^1(E,j_*L)\simeq W_T:=H_T/E_T,$ where
	\begin{multline*}
	H_T=\{(v_1,\ldots,v_{r+2})\in V^{r+2}\mid v_i\in \im(A_i-1),\\ 
	\quad v_1A_2\cdots A_r+\cdots +v_r+v_{r+1}(B-1)B^{-1}A^{-1}+v_{r+2}(1-A)B^{-1}A^{-1}=0\}	
	\end{multline*} 
	and
	\[
	E_T=\{(v(A_1-1),\ldots, v(A_r-1),v(A-1),v(B-1))\mid v\in V\}.
	\] 
\end{lem}

\subsection{Monodromy of Convolution Cocycles}\label{subsec:MonConv}
\begin{figure}[h!]
	\centering \tiny
	\begin{subfigure}{0.45\textwidth}
		\def\svgwidth{0.95\linewidth}
		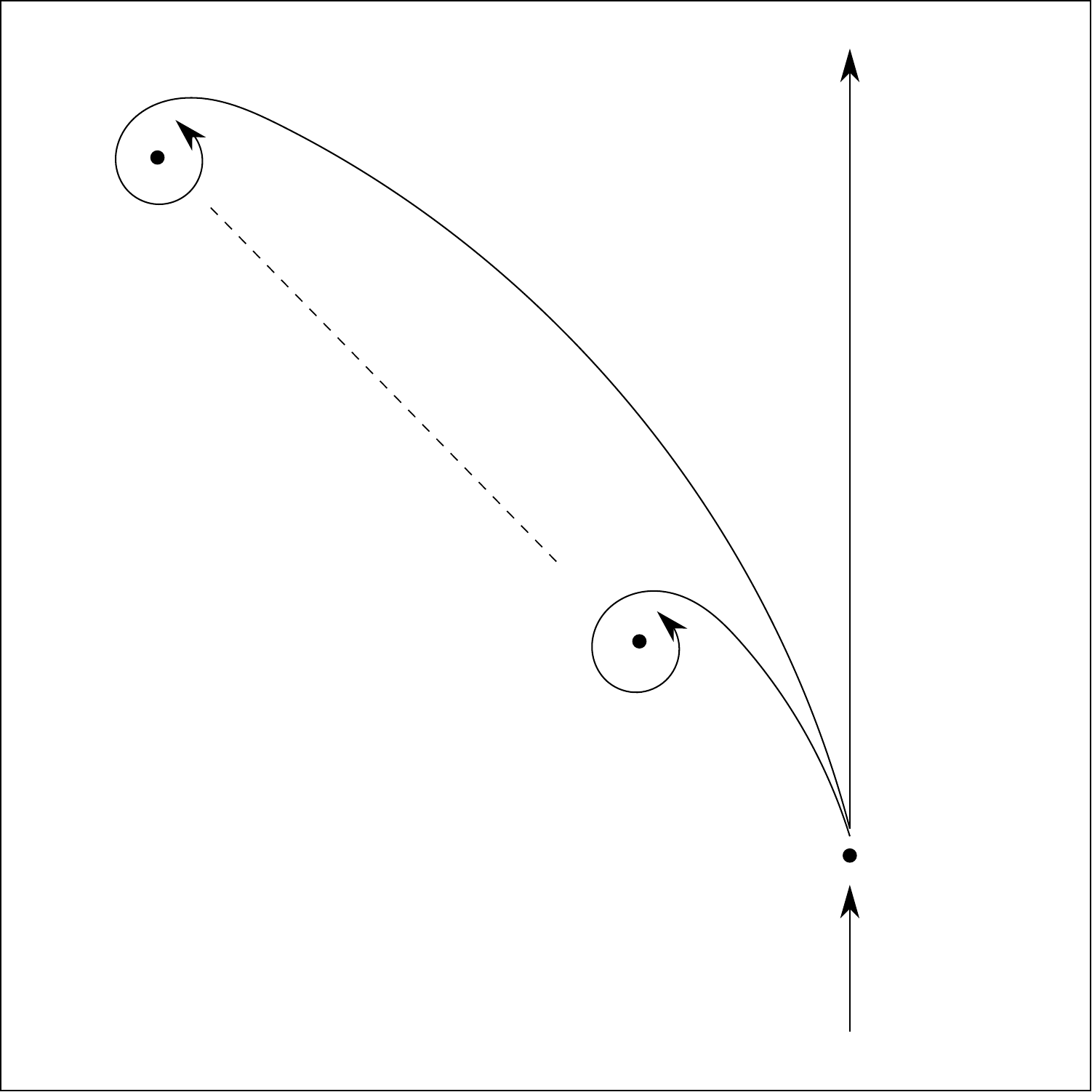{}
		\caption{Homotopy generators of the fibre}\label{fig:Braid_Conf2}
	\end{subfigure}
	\begin{subfigure}{0.45\textwidth}
		\def\svgwidth{0.95\linewidth}
		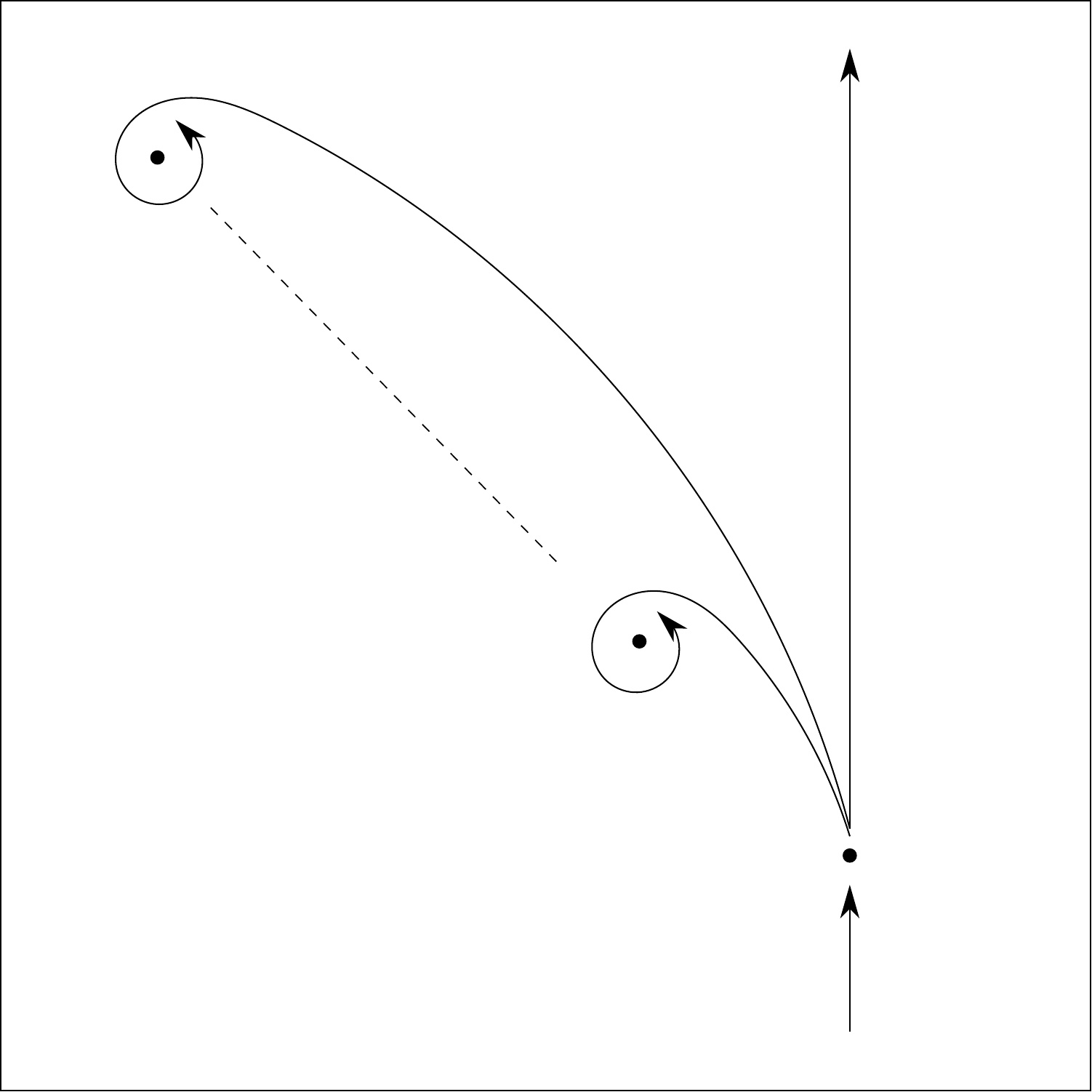{}
		\caption{Homotopy generators of the base.}\label{fig:Braid_Conf_Base}
	\end{subfigure}	
	\caption{The chosen setup} \label{fig:fig2}
\end{figure}

Let $\x_1:=\{x_1,\ldots,x_p\}\subseteq E$ and $\x_2:=\{y_1,\ldots,y_q\}\subseteq E$ and let throughout the rest of the article 
$L_i\,(i=1,2)$ be irreducible local systems on $U_i:=E\setminus \x_i,$ \emph{such that the associated 
	perverse intermediate extensions $K_i:=j_*L_i[1]\ (i=1,2)$ are in $\PPE$} (which is equivalent 
to the existence of at least one nontrivial local monodromy element $A_i$). Let
\begin{align}
(A_1,\ldots,A_p,A,B)\in \GL(V_1)^{p+2}\\
(B_1,\ldots,B_q,C,D)\in \GL(V_2)^{q+2}
\end{align}
be the respective monodromy tuple of $L_1$ and $L_2$ (with respect to 
monodromy generators as indicated in Fig.~\ref{fig:Braid_Conf2}). 

\medskip

We now compute as announced the monodromy representation $\eta$ of Eq.~\eqref{eq:ConvMonRep}:
\[ 
\eta:\pi_1(E\setminus \x_1\ast \x_2)\to \GL(H^1(E,j_*(L_1\otimes L_2(y_0-x))))
\]
by computing the monodromy tuple of the generic fibre of $L_1\otimes L_2(t-x)$ via
the braid identification $\phi_*\colon\pi_1(E\setminus \x_1\ast \x_2)\to P(p+q,E)$ obtained in \S \ref{secellbraids}. The homotopy generators of the fibre and the base are given as in Fig.~\ref{fig:fig2}. 

In order to do so, it is necessary to compute the monodromy tuple of $L_2(y_0-x),$ which is the same as the monodromy tuple of $L_2(-x)$ by neglecting the shift by $y_0.$
\begin{lem}\label{lemtupel} Let $L$ be a local system on $E\setminus \x\, (\x=\{x_1,\ldots,x_r\})$ having monodromy tuple $T$ according to the monodromy generators 
	$\gamma_1,\ldots,\gamma_r,\alpha,\beta$ chosen as in Fig.~\ref{fig:Braid_Conf2} with $\alpha_i$ replaced by $\gamma_i.$ 
	Let $\gamma_1',\ldots,\gamma_r',\alpha',\beta'$ be 
	similar monodromy generators of $E\setminus -\x\, (-\x=\{-x_1,\ldots,-x_r\}).$ Let $\varphi:\pi_1(E\setminus \x,x_0)\to \pi_1(E\setminus -\x,-x_0)$ denote 
	the homomorphism induced by $-x:E\to E, x\mapsto -x$ (see Fig.~\ref{fig:Braid_Conf_Phi2}).  
	Then
	\begin{multline*}
	\Bigl(\gamma_1',\ldots,\gamma_r',\alpha',\beta'\Bigr)=\Bigl(\varphi(\gamma_r),\varphi(\gamma_{r-1})^{\varphi(\gamma_r)},\ldots,\varphi(\gamma_1)^{\varphi(\gamma_2)\cdots 
		\varphi(\gamma_r)},\\
	\varphi(\gamma_r)^{-1}\cdots \varphi(\gamma_1)^{-1}\varphi(\alpha)^{-1}, \varphi(\beta)^{-1}\varphi(\gamma_1)\cdots\varphi(\gamma_r)\Bigr).
	\end{multline*}
	Especially, the monodromy tuple $T_{L_1\otimes L_2(y_0-x)}\in \GL(V_1\otimes V_2)^{p+q+2}$ of $L_1\otimes L_2(y_0-x)$ with respect to monodromy generators $\alpha_1,\ldots,\alpha_n,\alpha,\beta $ as above
	is given by 
	$$ (A_1\otimes 1_{V_2},\ldots,A_p\otimes 1_{V_2},1_{V_1}\otimes B_q, 1_{V_1}\otimes B_{q-1}^{B_q},\ldots, 1_{V_1}\otimes B_1^{B_2\cdots B_q},
	A\otimes \tilde{C} ,B\otimes \tilde{D}),$$
	where $1_{V_i}\, (i=1,2)$ denotes the identity  in $\GL(V_i),$  where  the $\otimes$-sign 
	denotes the usual Kronecker product of matrices and where
	$$ \tilde{C}:= B_q^{-1}\cdots B_1^{-1}C^{-1}\quad \text{ and } \quad \tilde{D}:=
	D^{-1}B_1\cdots B_q.$$ 
\end{lem}

\begin{proof} The first claim follows from the effect of $-x$ on the chosen homotopy generators and their translation back into our standard setup, as indicated by Fig.~\ref{fig:ConvG0}
	(note that if we represent $E$ by a square as above such that the focal point represents the origin of $E,$ then the $-x$-map is represented by the reflection 
	at the focal point, leading to Fig.~\ref{fig:Braid_Conf_Phi}).
	\begin{figure}[h!]
		\centering \tiny
		\begin{subfigure}{0.45\textwidth}
			\def\svgwidth{0.95\linewidth}
			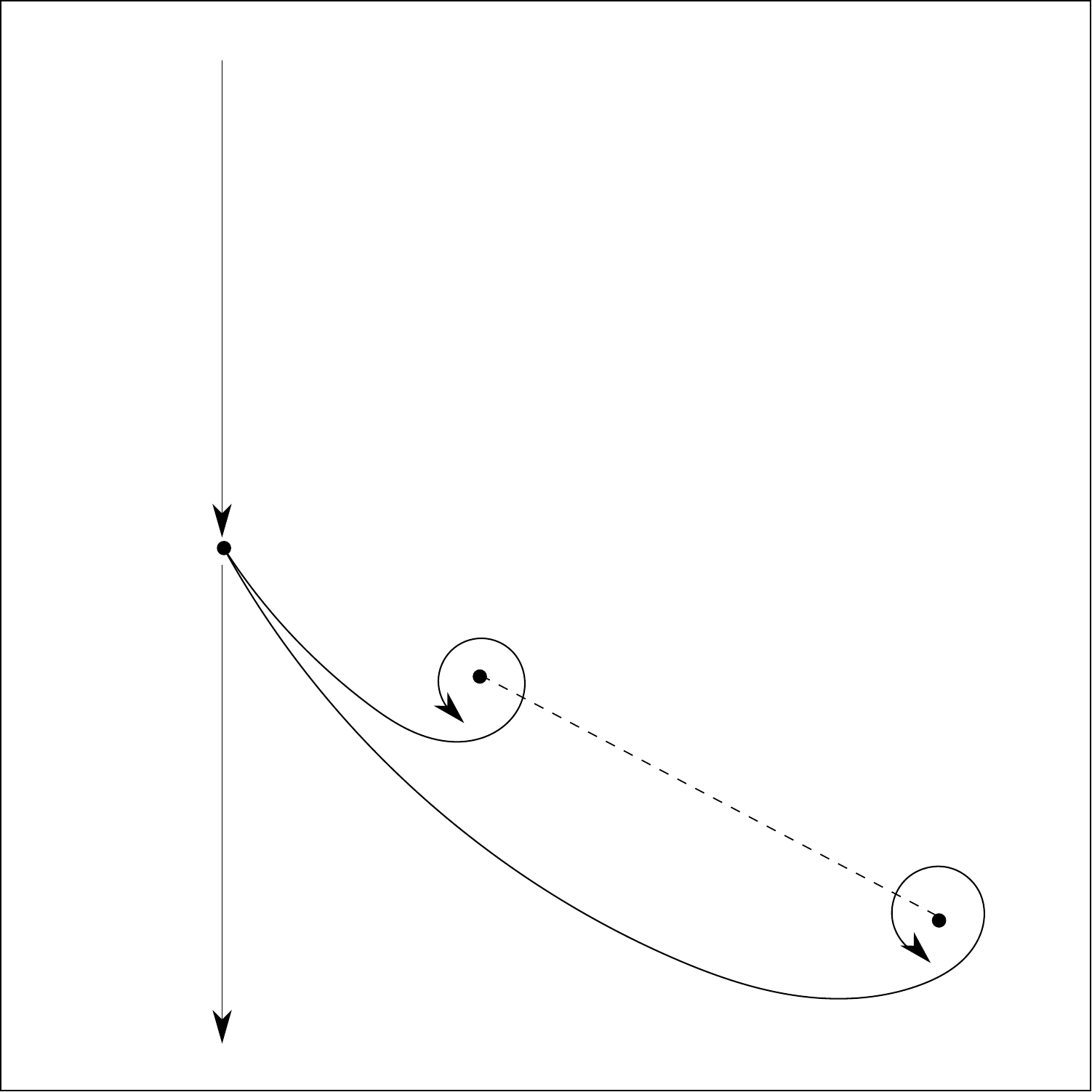{}
			\caption{Braids on $E\setminus\x$.}\label{fig:Braid_Conf_Phi}
		\end{subfigure}
		\begin{subfigure}{0.45\textwidth}
			\def\svgwidth{0.95\linewidth}
			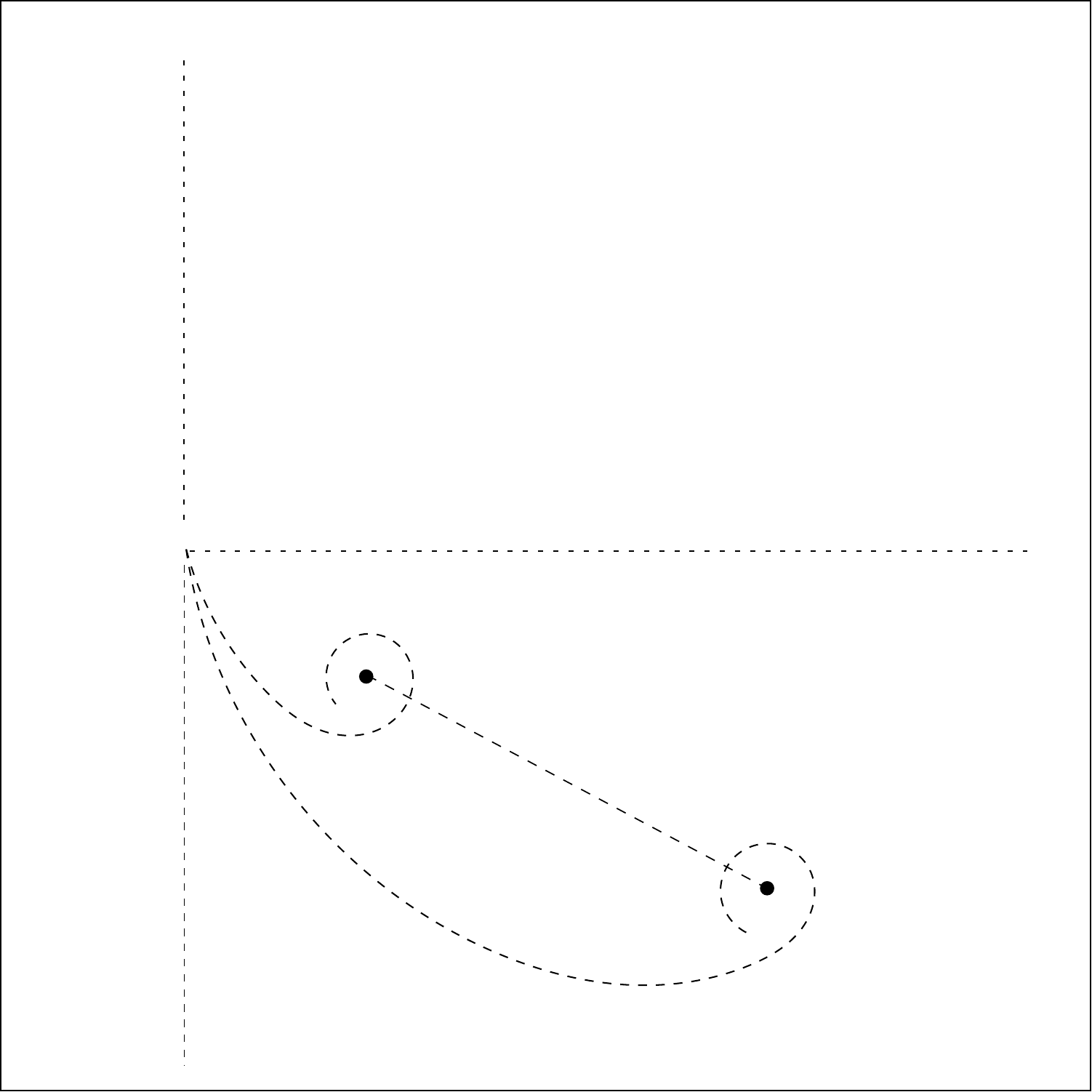{}
			\caption{Braids on $E\setminus-\x$}\label{fig:Braid_Conf_Phi2}
		\end{subfigure}	
		\caption{Effect of $-x$}
	\end{figure}
	The second claim is a direct consequence of the first one.
\end{proof}

One deduces the monodromy of the deformation of $\tau_{L_1\otimes L_2(y_0-x)}$ under the operation of $ \phi(\hat{\alpha})$ and  $\phi(\hat{\beta})$ as given by Prop.~\ref{propdeltt2}.

\begin{cor}\label{coreg} Let $\tau_{L_1\otimes L_2(y_0-x)}^{\phi(\hat{\alpha})}$ and 
	$\tau_{L_1\otimes L_2(y_0-x)}^{\phi(\hat{\beta})}$ be the deformations of $\tau_{L_1\otimes L_2(y_0-x)}$ under 
	$ \phi(\hat{\alpha}),  \phi(\hat{\beta})$.
	Then the componentwise application of the monodromy representation $ \rho_{L_1\otimes L_2(y_0-x)}$ of $L_1\otimes L_2(y_0-x)$ to $\tau_{L_1\otimes L_2(y_0-x)}^{\phi(\hat{\alpha})}$ and 
	$\tau_{L_1\otimes L_2(y_0-x)}^{\phi(\hat{\beta})}$ is given by (resp.)
	\begin{equation*}
	\begin{split}  
	\Bigl(A_1\otimes 1_{V_2},\ldots,A_p\otimes 1_{V_2},1_{V_1}\otimes B_q ,1_{V_1}\otimes B_{q-1}^{B_q},\ldots, 1_{V_1}\otimes B_1^{B_2\cdots B_q},
	A\otimes \tilde{C},B\otimes
	\tilde{D}\Bigr)^{1_{V_1}\otimes \tilde{C}^{-1}},\\
	\Bigl(A_1\otimes 1_{V_2},\ldots,A_p\otimes 1_{V_2},1_{V_1}\otimes B_q, 1_{V_1}\otimes B_{q-1}^{B_q},\ldots, 1_{V_1}\otimes B_1^{B_2\cdots B_q},A\otimes \tilde{C},B\otimes\tilde{C}\Bigr)^{1_{V_1}\otimes \tilde{D}^{-1}},
	\end{split}
	\end{equation*}
	using componentwise conjugation. 
\end{cor}

\begin{proof}  By Prop.~\ref{propdeltt2}, one has 
	$$ (\alpha_{1}^{\phi(\hat{\alpha})},\ldots, \alpha_{p+q}^{\phi(\hat{\alpha})} ,\alpha^{\phi(\hat{\alpha})},\beta^{\phi(\hat{\alpha})})=
	(\alpha_1^{\alpha_{p+1}\cdots\alpha_{p+q}},\ldots, \alpha_p^{\alpha_{p+1}\cdots\alpha_{p+q}}, \alpha_{p+1}^\alpha,\ldots, \alpha_{p+q}^\alpha,\alpha, \alpha_{p+q}^{-1}\cdots \alpha_{p+1}^{-1}\beta).$$ 
	Hence $ \rho_{L_1\otimes L_2(y_0-x)}(\alpha_{1}^{\phi(\hat{\alpha})},\ldots, \alpha_{p+q}^{\phi(\hat{\alpha})} ,\alpha^{\phi(\hat{\alpha})},\beta^{\phi(\hat{\alpha})})$ 
	is equal to 
	$$
	((A_1\otimes 1_{V_2})^{1_{V_1}\otimes (B_1\cdots B_q)},\ldots,(A_p\otimes 1_{V_2})^{1_{V_1}\otimes (B_1\cdots B_q)}, 
	(1_{V_1}\otimes B_q)^{A\otimes \tilde{C}},\ldots \quad \quad\quad \quad\quad \quad \quad \quad\quad \quad\quad \quad$$
	$$ \quad \quad\quad \quad\quad \quad \quad \quad\quad \quad\quad \quad \ldots  ,  (1_{V_1}\otimes B_1^{B_2\cdots B_q })^{A\otimes \tilde{C}}, A\otimes \tilde{C}, (1_{V_1}\otimes (B_1\cdots B_q)^{-1})\cdot
	(B\otimes \tilde{D}))=$$ 
	$$
	((A_1\otimes 1_{V_2})^{1_{V_1}\otimes \tilde{C}^{-1}},\ldots,(A_p\otimes 1_{V_2})^{1_{V_1}\otimes \tilde{C}^{-1}}, 
	(1_{V_1}\otimes B_q)^{1_{V_1}\otimes \tilde{C}^{-1}},\ldots \quad \quad\quad \quad\quad \quad \quad \quad\quad \quad\quad \quad$$ $$ \quad \quad\quad \quad\quad \quad \quad \quad\quad \quad\quad \quad \ldots  ,  (1_{V_1}\otimes B_1^{B_2\cdots B_q })^{1_{V_1}\otimes \tilde{C}^{-1}}, (A\otimes \tilde{C})^{1_{V_1}\otimes \tilde{C}^{-1}}, 
	(B\otimes \tilde{D})^{1_{V_1}\otimes \tilde{C}^{-1}}),$$
	where for the equality of the last component we argue as follows (the equalities for the other components hold trivially): it suffices to show that 
	$$ B_q^{-1}\cdots B_1^{-1} \tilde{D}=\tilde{C}\tilde{D}\tilde{C}^{-1},$$ 
	but this holds by the elliptic product relation. For $\tau_{L_1\otimes L_2(y_0-x)}^{\phi(\hat{\beta})}$ we use an analogous argument.
\end{proof}

The following result gives a complete description of the action on the fiber in terms of the cocyle algebra $W_{\bullet}=H_\bullet/E_\bullet$ and of the monodromy matrices $\{A_i\}$, $\{B_j\}$ and $A,\ B,\ C,\ D$ as defined in the Corollary above.

\begin{thm}  \label{etalem} Let $S:=E\setminus \x_1\ast \x_2$ 
	and let $U_0:= E\setminus (\x_1\cup y_0-\x_2).$
	Let $\gamma\in\pi_1(S)$ and $\delta:\pi_1(U_0)\to V_1\otimes V_2$ be a parabolic
	cocycle inside $H^1(U_0,L_1\otimes L_2(y_0-x)).$ We write $[\delta]$ for the class of $\delta$ in $W$. 
	Also, we choose a lift ${\tilde{\gamma}}\in\pi_1(\mathbb{V})$ of $\gamma$ according to the splitting in the diagram~\eqref{splitdiag}. Let 
	\[
	\eta:\pi_1(S)\to \GL(H^1(E,j_*(L_1\otimes L_2(y_0-x))))
	\]
	be the monodromy representation of 
	$L_1\ast L_2$ as above.
	Then
	$[\,\delta\,]^{\eta(\gamma)}=[\,\delta'\,]$, where
	$\delta':\pi_1(U_0)\to V$ is the cocycle
	\[
	\alpha \;\longmapsto\;  \delta(\tilde{\gamma}\alpha{\tilde{\gamma}}^{-1})\cdot \rho_{L_1\boxtimes L_2}(\tilde{\gamma})=\delta(\phi(\gamma)\alpha{\phi(\gamma)}^{-1})\cdot \rho_{L_1\boxtimes L_2}(\tilde{\gamma}), \qquad
	\alpha\in\pi_1(U_0). 
	\]
	Especially, a matrix $\Phi(\gamma)\in \GL((V_1\otimes V_2)^{p+q+2})$ which represents $\eta(\gamma)\in \GL(W_{T_{L_1\otimes L_2(y_0-x)}})$ is given by 
	the linear transformation, which transforms an evaluated cocycle 
	$$(\delta(\alpha_1),\ldots,\delta(\alpha_{p+q}),\delta(\alpha),\delta(\beta))$$
	into 
	$$(\delta(\alpha_1^{\phi(\gamma)^{-1}}),\ldots,\delta(\alpha_{p+q}^{\phi(\gamma)^{-1}}),\delta(\alpha^{\phi(\gamma)^{-1}}),\delta(\beta^{\phi(\gamma)^{-1}}))\cdot \rho(\tilde{\gamma})$$
	under the cocyle rule~\eqref{cocyclerule}, where the operation of $\gamma$  on 
	$\alpha_1,\ldots,\alpha_{p+q},\alpha,\beta$ is determined by Prop.~\ref{propdeltt} and Prop.~\ref{propdeltt2}, and where 
	$$ \rho(\tilde{\gamma})=\left\{\begin{array}{l} 1_{V_1}\otimes 1_{V_2} \quad \textrm{for $\gamma=\delta_{i,j}$}\\
	{1_{V_1}\otimes \tilde{C}^{-1}}  \quad \textrm{for $\gamma=\hat{\alpha}$}\\
	{1_{V_1}\otimes \tilde{D}^{-1}}  \quad \textrm{for $\gamma=\hat{\beta}$}\end{array} \right. , $$
	acts componentwise. 
\end{thm}

\begin{proof} The first claim is completely analogous to \cite{dw}, Lem.~2.2. The expression of $ \rho(\tilde{\gamma})$ 
	follows from the definition of $L_1\otimes L_2(y-x)$ by projecting $\tilde{\gamma}$ via the $y-x$-map to $E\setminus \x_2.$	\end{proof}

\section{Seven-point Sheaves with Tannakian group equal to $G_2$} \label{secbeauville}
We establish the main result of this paper, which states that the exceptional algebraic group $G_2\leqslant \GL_7$ is recovered as the Tannaka group of a certain class of seven-point sheaves in $\PPE$ on elliptic curves $E$ (see Thm.~\ref{thmapp}).

\subsection{Seven-point Sheaves over Elliptic Curves}\label{sevenpoint}

The following definition generalizes the geometric construction obtained via the Beauville classification of elliptic surfaces with four bad semistable fibres \cite{Beauville,KatzEll}.

\begin{defn}\label{def:7shef}
	Let $E$ be an elliptic curve over an algebraically closed field $\bar k$, and let $\ell\neq \char(\bar k)$. A \emph{seven-point sheaf} $N$ on $E$ is an irreducible $\bQl$-perverse sheaf $N\in \PPE$ which is self-dual, and smooth over $E(k)$ except at seven points of the form $\mathbf{a}=\{\pm a_1,\pm a_2,\pm( a_1+a_2),0\}$ over which it has local monodromy of the form:
	\begin{equation}\label{eq:mon7}
	(A_4,A_3^{A_4},A_2^{A_3A_4},A_1^2,A_2,A_3,A_4)\quad A_1\cdots A_4=1 
	\end{equation}
	where $A_i\in SL_2$ $(i=1,...,4)$ are nontrivial tame and unipotent, and with trivial handle monodromy.
\end{defn}

Following the proof of Lem.~\ref{lem:7sheafP1}, one shows that the $A_i$ ($i=1,\dots, 4)$ are given as in Tab. \ref{tab:locMon}.

\begin{table}[h!]\centering
	\begin{tabular}{p{.5em}p{1em}cccp{5em}}\toprule
		& &\multicolumn{3}{c}{Local monodromy}& Condition\\  
		& & $A_1$ & $A_2$ & $A_3$ & \\ \toprule 
		I & & $\begin{pmatrix}
		1& 1  \\
		0   &1
		\end{pmatrix}$ &
		$\begin{pmatrix}
		1& -1  \\
		0   &1
		\end{pmatrix}$ &
		$\begin{pmatrix}
		1& 0  \\
		y  &1
		\end{pmatrix} $ & $y\neq 0$\\ \midrule 
		& i)& $\begin{pmatrix}
		1& 1  \\
		0   &1
		\end{pmatrix}$ &
		$\begin{pmatrix}
		1& 0  \\
		\frac{(a-1)^2}{b(a+b)}  &1
		\end{pmatrix}$ &
		$\begin{pmatrix}
		a& b  \\
		-\frac{(a-1)^2}{b}  &\,\,2-a
		\end{pmatrix}$ & $a\neq 1$, $b\neq 0$, $a+b\neq 0$ \\ \cmidrule(l){2-6}
		II & ii)& $\begin{pmatrix}
		1& 1  \\
		0   &1
		\end{pmatrix}$ &
		$\begin{pmatrix}
		1& 0  \\
		y  &1
		\end{pmatrix}$ &
		$\begin{pmatrix}
		1& 0  \\
		-y &1
		\end{pmatrix}$ & $y\neq 0$\\ \cmidrule(l){2-6}
		& iii) & $\begin{pmatrix}
		1& 1  \\
		0   &1
		\end{pmatrix}$ &
		$\begin{pmatrix}
		1& 0  \\
		y  &1
		\end{pmatrix}$ &
		$\begin{pmatrix}
		1& -1  \\
		0 &1
		\end{pmatrix}$ &$y\neq 0$\\
		\bottomrule	 		   
	\end{tabular}
	\caption{Seven-point sheaves local monodromies}\label{tab:locMon}
\end{table}

An example of seven-point sheaves is given by the construction in the following lemma, that is similar to Katz's one via Beauville classification.

\begin{lem}\label{lem:7sheafP1}
	Let $L$ be an irreducible constructible $\ell$-adic $\bQl$-sheaf on $\PP^1_\CC$ which is smooth on $\PP^1_\CC\setminus \{x_1,x_2,x_3,x_4=\infty\}$, with nontrivial tame unipotent local monodromy at each $x_i\, (i=1,\ldots,4)$, and such that the singular fibre of the elliptic involution quotient $\pi\colon E\to \PP^1$ has cardinality $7$ of the form:
	\[
	\pi^{-1}(\{x_1,x_2,x_3,x_4\})=\{\pm a_1,\pm a_2,\pm( a_1+a_2),0\}\subseteq E(\CC).
	\]		
	Then $N:=\pi^*L[1]\in \PPE$ is a seven-point sheaf. 
\end{lem}

Reciprocally, a seven-point sheaf $N$ on $E$ descends to a $\ell$-adic sheaf $L$ on $\PP^1$ as above.

\begin{proof}
	Let $L$ on $\PP^1(\CC)\setminus \{x_1,\ldots,x_4\}$ be as above, and let us choose an isomorphism 
	$\bQl\simeq \CC$. One obtains that the set of such irreducible $\bQl$-local systems with unipotent local monodromy at $x_i (i=1,\ldots, 4\}$ is parametrized as follows: After fixing generators $\gamma_1,\ldots,\gamma_4$ of $\pi_1(\PP^1(\CC)\setminus \{x_1,\ldots,x_4\})$ satisfying $\gamma_1\cdots\gamma_4=1,$
	any such local system is determined by its monodromy tuple $(A_1,\ldots,A_4)\in \SL_2(\CC)^4,$ satisfying $A_1\cdots A_4=1.$  By conjugation, one may assume 
	$A_1$ to be the standard unipotent upper triangular matrix such that the centralizer $C_{\SL_2(\CC)}(A_1)$ consists of the unipotent upper triangular matrices, up to 
	$\pm1.$ By conjugating the second element with a suitable element in 
	$C_{\SL_2(\CC)}(A_1)$ and using the trace-2-condition for $A_1A_2A_3$ we obtain that we are in one of the cases above.
	
	One checks immediately that the quadratic pullback gives a perverse sheaf $N=\pi^*L[1]\in \PPE$ over $E$ with monodromy tuple as in Eq.~\eqref{eq:mon7}, thus a seven-point sheaf over $E$.
\end{proof}

\subsection{$G_2$ as Tannaka group over $E$}
The characterization of the Tannaka group $G_N$ of a seven-point sheaf $N\in \PPE$ relies on $N\ast N$-decomposition and Tannakian properties.

\begin{prop}\label{propm3} For any seven-point $\bar{\QQ}_{\ell}$-sheaf $N\in \PPE$ on $E_{\bar k}$ over an algebraically closed field $\bar k$, the convolution $N* N$ contains a summand isomorphic to $N$. 
\end{prop}

\begin{proof} Assume that $\bar k=\CC$ and choose an isomorphism $\bQl\simeq \CC$. We first prove this result for analytic local systems $N^\an\in \PPE$ over $E(\CC)\setminus \{\pm a_1,\pm a_2,\pm( a_1+a_2),0\}$, by describing the monodromy tuple of $N^\an\ast N^\an$ in $\pi_1(E\setminus \x_1\ast \x_2,y_0)$ as in \S~\ref{subsec:MonConv}, where we write $\x_1:=\{x_1=-(a_1+a_2),x_2=-a_2,x_3=-a_1,x_4=0,x_5=a_1,x_6=a_2,x_7=a_1+a_2\}$ and $\x_2:=\{y_1:=-x_1,\ldots,y_7:=-x_7\}.$

	By Rem.~\ref{Rem:genVSdeg} the computation can be carried out in the standard $G_2$-configuration with
	$$\#\x_1\ast \x_2=\#\{x_i+x_j\mid i,j=1,\ldots ,7\}=19$$
	and a singularity arrangement of $N\boxtimes N$ of the form indicated in Fig.~\ref{fig:ConvolEll}.
	\begin{figure}[h!]
		\centering \small
		\def\svgwidth{0.95\linewidth}
		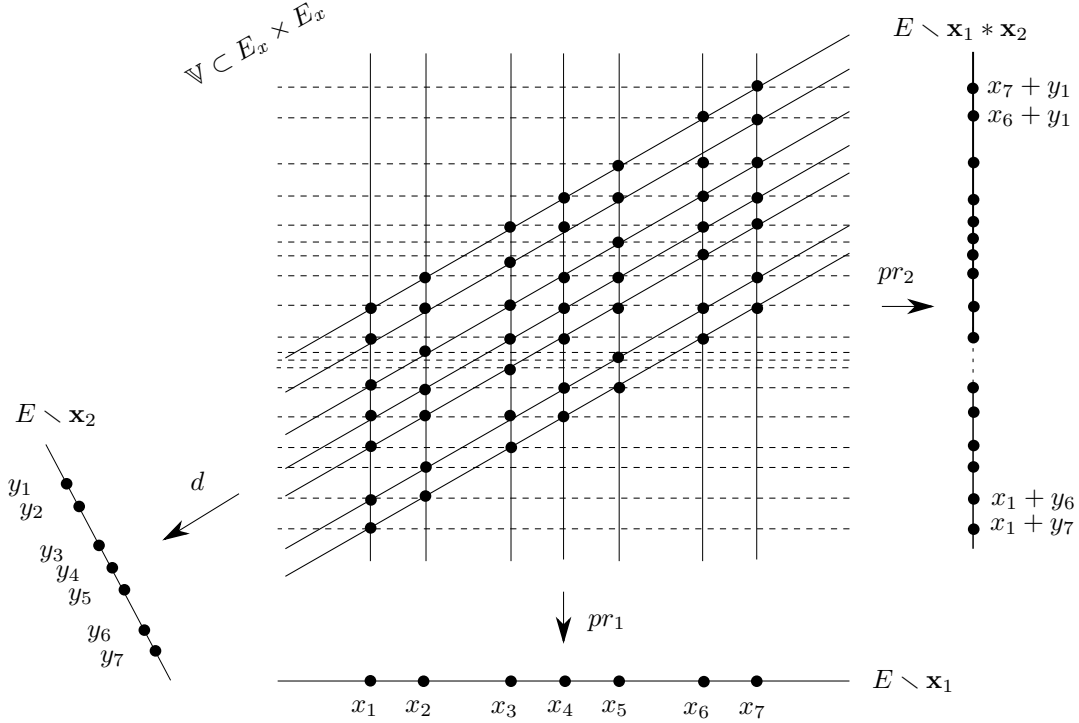{}
		\caption{Convolution of a seven-point sheaf over $E$.}\label{fig:ConvolEll}
	\end{figure}
	Using this one obtains the expression of the paths $\delta_1=\delta_{1,7},\delta_2=\delta_{2,7},\ldots, \delta_{18}=\delta_{6,1},\delta_{19}=\delta_{7,1}\in \pi_1(E_\setminus \x_1\ast \x_2,y_0)$ 
	in terms of the braid group generators $\beta_1,\ldots,\beta_{13}$ (with the conventions of \S~\ref{secellbraids}) as presented in Appendix~\ref{app:braids}.
	
	Applying Formula~\eqref{eqbetai} one obtains the associated deformations of the standard generators $\alpha_1,\ldots,\alpha_{14}.$ Similarly, one uses  
	Prop.~\ref{propdeltt2} in order to get the corresponding deformations under the global braids $\hat{\alpha}$ and $\hat{\beta}$ of $E\setminus \x_1\ast \x_2.$ 
	The last claim of Thm.~\ref{etalem} and an implementation of this result into MAGMA-coding   \cite{MAGMA}, also presented in Appendix \ref{app:MonMagma},
	gives the monodromy matrices of $N^\an\ast N^\an.$ Decomposing the associated representation, one shows that in all cases one has a decomposition
	\begin{equation}\label{eq:decomNN}
	N^\an \ast N^\an= N_1\oplus N_2\oplus N_3\oplus \delta_0,
	\end{equation}
	where $N_2,N_3$ are irreducible of generic rank $8$ and $18$ (resp.), and $N_1$ is isomorphic to $N ^\an$ as announced.

	Because of the equivalence of categories between $\overline{\mathbb Z}_\ell$-local systems in the \'{e}tale and analytic topologies, which commutes with cohomology, the same property holds for $\ell$-adic local systems on $E(\CC)$. 
	
	\medskip

	We now show that the same property holds over any arbitrary algebraically closed base field $\bar{k}$. First we handle the case of $char(\bar k)=0$. The \'{e}tale cover in question has monodromy contained in a compact subgroup of $GL_2( \overline{\mathbb Z_\ell} )$, which is necessarily a countable inverse limit of finite groups, because any compact subgroup of $GL_2 ( \overline{\mathbb Z_\ell} / \ell^n)$ is finite. Hence it can be defined using a countable system of finite \'{e}tale covers, which are each defined over a subfield of $\bar k$ finitely generated over $\mathbb Q$, hence all together defined over a subfield of $\bar k$ of countable transcendence degree, which necessarily embeds into $\CC$. Because base-changing between this subfield and $\bar k$ -- or this subfield and $\CC$ -- preserves cohomology groups, $N\ast N$ has the same summands over an arbitrary algebraically closed field of characteristic zero.
	
	\medskip
	
	Next we handle the case of $char(\bar k)=p>0$. Considering $R$ the Witt vectors ring of $\bar k$, we lift $E_{\bar k}$ arbitrarily to an elliptic curve over $R$ and $x_1,x_2,x_3,x_4=\infty$ to $R$-points of $\mathbb P^1_R$. We can lift the local system $L:=R^1\pi_* N$ on $\PP_{\bar k}\setminus\{x_1,x_2,x_3,x_4\}$ to a local system on $\mathbb P^1_R \setminus\{x_1,x_2,x_3,x_4\}$, because it has tame local monodromy, and the tame fundamental group in characteristic $p$ is a quotient of the tame fundamental group in characteristic zero.  Hence we obtain $N$ as a complex of sheaves on $E_R$, that is perverse on the generic and special fibers and has tame local monodromy around its singular points $\{\pm a_1, \pm a_2, \pm (a_1+a_2),0\}$. The convolution $N * N$ may be defined relative to $R$ as the pushforward of $N \boxtimes N$ along the smooth proper  multipliation map $(E \times E)_R \to E_R$.  The sheaf $N \boxtimes N$ is lisse away from the divisor 
	\[ \pm a_1 \times E\cup \pm a_2 \times E\cup \pm (a_1+a_2) \times E \cup 0 \times E \cup E \times \pm a_1\cup E \times \pm a_2\cup E \times \pm (a_1+ a_2) \cup E \times 0 .\] 
	It has tame ramification around that divisor. Hence by Deligne's semicontinuity Theorem in \cite{LAU81} Cor.~2.1.2, $N \ast N $ is lisse on the open set where this divisor is smooth and transverse to the multiplication map. Each irreducible component is smooth and transverse, so this open set is exactly the set where the irreducible components do not intersect, i.e. the complement of the sums $\{ y_1+y_2| y_1, y_2 \in \{\pm a_1, \pm a_2, \pm (a_1+a_2),0\}\}$. (The semicontinuity theorem covers the case of the extension by zero of such a lisse sheaf, while $N \boxtimes N$ is nonzero on this divisor. However, it is lisse on the divisor, so we can write it as the extension of a lisse sheaf supported on that divisor by an extension by zero of a lisse sheaf on its complement, showing that its pushforward is lisse).
	
	The map $N \to N \ast N$ is a map between sheaves lisse on the open complement inside $E_R$ of $\{ y_1+y_2| y_1, y_2 \in \{\pm a_1, \pm a_2, \pm (a_1+a_2),0\}\}$. Previously, we showed this map existed and was injective on the geometric generic fiber of this open set. Because both sheaves are lisse, it extends to the whole open set, and remains injective on the special fiber. Because $N$ is irreducible on the special fiber, it remains injective, showing that $N$ is a summand of $N\ast N$ on the special fiber as well, as desired.
	
\end{proof}

\begin{rem}\label{Rem:genVSdeg} Let 
	$$\widehat{G}_2=\{\mathbf{a}=\{\pm a_1,\pm a_2,\pm(a_1+a_2),0\}\mid a_1,a_2\in E(\CC),\, \#\mathbf{a}=7 \}$$ be the configuration space belonging to the group $G_2$. It has a Zariski dense subset of points $\x$ such that $\#\x\ast \x=19$ (the generic case). Let $\mathbf{a}_0 \in \widehat{G}_2$ be a generic and real configuration as indicated in Fig.~\ref{fig:ConvolEll}. If $\mathbf{b}$ is any point of $\widehat{G}_2$, it can be connected by a path $\sigma\colon [0,1]\to \widehat{G}_2$ starting from $\mathbf{a}_0$ and ending on $\mathbf{b}$, such that $\sigma(t)$ is generic $\forall t\in [0,1)$. This path can be given as the effect of a family of homeomorphisms $I_t\colon E\to E$ such that: $I_0=id$, $\sigma(t)=I_t(\mathbf{a}_0)$, and $I_t$ gives the identity on some specific area (depending on $\mathbf{b}$ and not on $t$). 
	Let $N_t$ denote the pullback of a given $7$-point sheaf $N_0$ with singular locus $\mathbf{a}_0$ along the homeomorphism $I_t^{-1}.$ 
	The family of homeomorphisms $I_t$ induces a family of homeomorphisms $$\mathcal{I}_t\colon E\times E\setminus \{\text{convol.  arrangement wrt } \mathbf{a}_0 \}\to E\times E\setminus\{\text{convol.  arrangement  wrt } I_t(\mathbf{a}_0)\},$$ 
	where the convolution arrangement is defined with respect to the self convolution of $N_t$ and is the inverse image of the singular loci under the two projections and the substraction map, cf.~Fig.~\ref{fig:ConvG0}. One obtains the same monodromy matrices for $N_t\ast N_t$, $t\in [0,1),$
	if we carry the whole initial topological setup (homotopy generators of the fibre, braid group generators, etc.) along, using the family of homeomorphisms $\mathcal{I}_t.$ 
	
	In the degenerate case, some of the monodromy matrices are multiplied up according to the collapsing of points on $\mathbf{b}\ast\mathbf{b}$, so that the monodromy group in the degenerate case is a subgroup of the monodromy group of the non-degenerate case, and we conclude as in Prop.~\ref{g2criterion}.
	
	Such degeneracies are given by $\#\mathbf{b}\ast \mathbf{b}=n$ with $n\in\{7,8,9,\dots, 18\}$ and include cases coming from torsion identities, i.e. with $a_1,a_2\in E(\CC)_{tors}$: for example $n=7$ is obtained with $a_2=2a_1$ and $a_1\in  E(\CC)[7]$.
\end{rem}

\begin{prop}\label{g2criterion} Let $K\in \PPE$ be an irreducible perverse sheaf such that: $\chi(K) =7$,  $\rk(K)$ is not divisible by $14$, $DK = [-1]^* K$, $K \ast  K$ contains a summand isomorphic to $K$, and $K$ is not isomorphic to its own translate by any nontrivial point of $E$. Then $K$ has Tannakian monodromy group $G_K$ equal to $G_2$. 
\end{prop}

\begin{proof} By Lem.~\ref{WS3}, $K$ corresponds to a seven-dimensional representation of $G_K$ which is self-dual, whose third tensor power contains a trivial representation, and which is irreducible when restricted to the identity component. It follows from 
	the assumption $DK = [-1]^* K$ that $G_K\leqslant O_7$. The existence of the summand $K$ inside $K\ast K$ 
	implies that 
	$\delta_0$ is a summand of $K\ast K \ast K,$ forming a nontrivial third moment of $K$ (this implies the existence of a fixed line in $H^0(E,K^{*3})$), which implies that $G_K\leq G_2$ 
	as explained in \cite{KatzEll}, \S~3.  Hence, by the
	classification of irreducible subgroups of $\SO_7$ (cf.~\cite{KatzEll}, \S~3), the Tannakian group of $K$ is either $G_2$, acting by its standard representation, or $PGL_2\leq G_2$, acting by its seven-dimensional irreducible representation. 
	
	Let us show a contradiction in the second case: let $F$ be the perverse sheaf corresponding by general Tannakian formalism to the three-dimensional irreducible representation of $PGL_2$. Then by Lem.~\ref{WS1}, one has 
	\[\begin{cases*}
	\rk(F\ast  F ) = 6 \rk (F)\\
	\rk(F \ast  F \ast F) = 27 \rk(F)
	\end{cases*}\quad \text{ while }\quad  F \ast  F \ast  F + \delta_0 = K + 2 F \ast  F + F 
	\]
	by an identity in the representation ring of $SL_2$, cf. \cite{OnischikVinberg}, Reference Chapter, Tab. 5 - $A_1$, so $\rk(K) = 14 \rk(F)$. 
\end{proof}

\begin{thm}\label{thmapp} Let $E$ be an elliptic curve over $\bar k$ and $N\in \PPE$ a seven-point sheaf on $E$. The Tannaka group $G_N$ of $N$ is isomorphic to the exceptional simple group of type $G_2.$ 
\end{thm}

\begin{proof}
	We check all the conditions of Prop.~\ref{g2criterion}: The Euler characteristic of $N$ is $7$ by the Euler characteristic formula; the rank of $N$ is $2$ (which is not divisible by $14$); the sheaf $N$ is isomorphic to both its dual and its pullback on the inverse map, so these are isomorphic to each other.
	
	For the translation property, note that because it has seven singular points, it could only be isomorphic to its translate by a point of order seven. In which case, the seventh power of this isomorphism would be a scalar multiplication. Adjusting the isomorphism, we can assume that the scalar is one. Combined with the inverse map, this gives a group isomorphic to $D_7$ of automorphisms of the elliptic curve that act also on the sheaf. The only automorphisms that have fixed points are the reflections, which are isomorphic to the standard reflection, which acts as the identity at the stalk on its fixed points, so all these isomorphisms act as the identity on their fixed points. Hence the group action gives descent data to the quotient $E / D_7 \cong \mathbb P^1$, producing an irreducible sheaf of rank two on $\mathbb P^1$ whose only singular point is the image of the singular points of $N$, with unipotent local monodromy at that point, which is absurd. 
	
	Finally, the summand condition on $N\ast N$ is given by Prop.~\ref{propm3}.
\end{proof}

\begin{rem}\label{rem:Katz} In the special case of seven-point sheaves $N\in \PPE$ on $E_{\bar k}$ which are built from $\PP^1$ using Beauville surface, this result has been first obtained by Katz in \cite{KatzEll}, Thm.~4.1. using Frobenius traces computations. Katz raises the questions of finding \emph{``a conceptual explanation''} and of ``seeing'' \emph{``an alternating trilinear form on this piece of $H^2$''} -- see ibid. Rem.~5.2. The fact that our monodromy approach recovers Katz' result (see Appendix \ref{app:BeauvMon} for the monodromy matrices in this case), and the fact that $N\ast N$ contains a copy of $N$ (see Eq.~\eqref{eq:decomNN}) -- leading to a $\delta_0\in N\ast N\ast N$ which, under $H^0$, gives the awaited trilinear form on $\Lambda^3(H^0(E,N))$--, provide some answers to those two questions.
\end{rem}

\newpage
\appendix

\section{Appendix -- MAGMA Braids Computations}\label{app:braids}

The following list gives in MAGMA-code the braids mentioned in the proof of Prop.~\ref{propm3}. 

\lstset{
	columns=flexible,breaklines=true, breakatwhitespace=true, 
	basicstyle=\small\ttfamily,	
	keywordstyle=\color{black}\bfseries\underbar,
	identifierstyle=,
	stringstyle=\ttfamily,
	showstringspaces=false, showspaces=false, showtabs=false,
	breakindent=2em,
	frame=single]}
\lstset{commentstyle=\color{Grey}\normalfont,
	language=Java, morekeywords={function},
}

\begin{lstlisting}
//The tuples encode the braids in a way that a positive number i corresponds to the
// braids \beta_i and such that a negative number -i corresponds to \beta_i^{-1}
delta:=[];
delta[1]:=[7, 6, 8, 7, 5, 9, 4, 6, 8, 10 ,3 ,5 ,9 ,11, 4, 7, 10, 2, 6, 8, 12, 1, 3, 5, 7, 9, 11, 13, 2, 6, 8, 12, 4, 7, 10, 3, 5, 9, 11, 4, 6, 8, 10, 5, 9, 7, 6, 8, 7, 7, -8, -6, -7, -9, -5, -10, -8,-6, -4, -11, -9, -5, -3, -10, -7, -4, -12, -8, -6,-2, -13, -11, -9, -7, -5, -3, -1, -12, -8, -6, -2, -10, -7, -4, -11, -9, -5, -3,-10, -8, -6, -4, -9, -5, -7, -8, -6, -7];
delta[2]:=[7,6,8,7,5,9,4,6,8,10,3,5,9,11,4,7,10,2,6,8,12,1,3,5,7,9,11,13,2,6,8,12, 4,7,10,3,5,9, 11,4,6,8,10,5,9,7,6,6,8,8,-7, -9, -5, -10, -8,-6, -4, -11, -9, -5, -3, -10, -7, -4, -12, -8, -6, -2, -13, -11, -9, -7, -5, -3, -1, -12, -8, -6, -2, -10, -7, -4, -11, -9, -5, -3, -10, -8, -6, -4, -9, -5, -7, -8, -6, -7];
delta[3]:=[7,6,8,7,5,9,4,6,8,10,3,5,9,11,4,7,10,2,6,8,12,1,3,5,7,9,11,13,2,6,8,12, 4,7,10,3,5,9, 11,4,6,8,10,5,9, 7,7, -9, -5, -10, -8,-6, -4, -11, -9, -5, -3, -10, -7, -4, -12, -8, -6, -2, -13, -11, -9, -7, -5, -3, -1, -12, -8, -6, -2, -10, -7, -4, -11, -9, -5, -3, -10, -8, -6, -4, -9, -5, -7, -8, -6, -7 ];
delta[4]:=[7,6,8,7,5,9,4,6,8,10,3,5,9,11,4,7,10,2,6,8,12,1,3,5,7,9,11,13,2,6,8,12, 4,7,10,3,5,9,11,4,6,8,10, 5,5,9,9, -10, -8,-6, -4, -11, -9, -5, -3, -10, -7, -4, -12, -8, -6, -2, -13, -11, -9, -7, -5, -3, -1, -12, -8, -6, -2, -10, -7, -4, -11, -9, -5, -3, -10, -8, -6, -4, -9,-5, -7, -8, -6, -7];
delta[5]:=[7,6,8,7,5,9,4,6,8,10,3,5,9,11,4,7,10,2,6,8,12,1,3,5,7,9,11,13,2,6,8,12, 4,7,10,3,5,9,11, 4,4,6,6,8,8,10,10,-11, -9, -5, -3, -10, -7, -4, -12, -8, -6, -2, -13, -11, -9, -7, -5, -3, -1, -12, -8, -6, -2, -10, -7, -4, -11, -9, -5, -3, -10, -8, -6, -4, -9, -5, -7, -8,-6, -7 ];
delta[6]:=[7,6,8,7,5,9,4,6,8,10,3,5,9,11,4,7,10,2,6,8,12,1,3,5,7,9,11,13,2,6,8,12, 4,7,10,3,3,5,5, 9,9,11,11, -10, -7, -4, -12, -8, -6, -2, -13, -11, -9, -7, -5, -3, -1, -12, -8, -6, -2, -10, -7, -4, -11, -9, -5, -3, -10, -8, -6, -4, -9, -5, -7, -8, -6, -7];
delta[7]:=[7,6,8,7,5,9,4,6,8,10,3,5,9,11,4,7,10,2,6,8,12,1,3,5,7,9,11,13,2,6,8,12, 4,10,7,7, -10,-4, -12, -8, -6, -2, -13, -11, -9, -7, -5, -3, -1, -12, -8, -6, -2, -10, -7, -4, -11, -9, -5, -3, -10, -8, -6, -4, -9, -5, -7, -8, -6, -7 ];
delta[8]:=[7,6,8,7,5,9,4,6,8,10,3,5,9,11,4,7,10,2,6,8,12,1,3,5,7,9,11,13,2,6,8,12, 4,4,10,10, -12, -8, -6, -2, -13, -11, -9, -7, -5, -3, -1, -12, -8, -6, -2, -10, -7, -4, -11, -9, -5, -3, -10, -8, -6, -4, -9, -5, -7, -8, -6, -7 ];
delta[9]:=[7,6,8,7,5,9,4,6,8,10,3,5,9,11,4,7,10,2,6,8,12,1,3,5,7,9,11,13,2,2,6,6,8, 8,12,12,-13, -11, -9, -7, -5, -3, -1, -12, -8, -6, -2, -10, -7, -4, -11, -9, -5, -3,-10, -8, -6, -4, -9, -5, -7, -8, -6, -7 ];
delta[10]:=[7,6,8,7,5,9,4,6,8,10,3,5,9,11,4,7,10,2,6,8,12,1,1,3,3,5,5,7,7,9,9,11, 11,13,13, -12, -8, -6, -2, -10, -7, -4, -11, -9, -5, -3, -10, -8, -6, -4, -9, -5, -7,-8, -6, -7];
delta[11]:=[7,6,8,7,5,9,4,6,8,10,3,5,9,11,4,7,10,2,2,6,6,8,8,12,12, -10, -7, -4, -11, -9, -5, -3, -10, -8, -6, -4, -9, -5, -7, -8, -6, -7 ];
delta[12]:=[7,6,8,7,5,9,4,6,8,10,3,5,9,11,7,4,4,10,10, -7,-11, -9, -5, -3, -10, -8, -6, -4, -9, -5, -7, -8, -6, -7 ]; 
delta[13]:=[7,6,8,7,5,9,4,6,8,10,3,5,9,11,7,7, -11, -9, -5, -3, -10, -8, -6, -4, -9, -5, -7, -8, -6, -7 ];
delta[14]:=[7,6,8,7,5,9,4,6,8,10,3,3,5,5,9,9,11,11,-10, -8, -6, -4, -9, -5, -7, -8, -6, -7];
delta[15]:=[7,6,8,7,5,9,4,4,6,6,8,8,10,10,-9, -5, -7, -8, -6, -7 ];
delta[16]:=[7,6,8,7,5,5,9,9,-7, -8, -6, -7 ];
delta[17]:=[7,6,8,7,7,-8,-6,-7];
delta[18]:=[7,6,6,8,8,-7];
delta[19]:=[7,7];
\end{lstlisting}

\section{Appendix - MAGMA Monodromy of the Convolution}\label{app:MonMagma}

The following lists the MAGMA code which is used in the proof of Prop.~\ref{propm3}. 
\pagestyle{empty}
\begin{lstlisting}
Invert:=function(W)
//W tuple describing a word in the homotopy generators
//(gamma_1,...,gamma_r)=(\alpha_1,...,\alpha_{p+q},\alpha,\beta))
//where \gamma_k is encoded as [k,1] and \gamma_k^-1 is encoded as [k,-1] 
//The function returns the tuple corresponding to the inverse word.
local I,i;
I:=Reverse(W);
for i in [1..#W] do
I[i][2]:=-1*I[i][2];
end for;
return I;
end function;

Conjugate:=function(W1,W2)
//Gives the effect of conjugation on tuples describing 
//words in the  homotopy generators 
//(gamma_1,...,gamma_r)=(\alpha_1,...,\alpha_{p+q},\alpha,\beta)) as above 
local C;
C:=(Invert(W1) cat W2) cat W1;
return C;
end function;

ConjugateInverse:=function(T1,T2)
//Gives the conjugation with the inverse
local C;
C:=(T1 cat T2) cat Invert(T1);
return C;
end function;

Braid:=function(j,W)
//Describes the effect of the right tuple-braiding action on of the braid \beta_j
//(cf. Section 2.1)
//on W=[[[1,1]],...,[[r,1]]] (W=Tuple describing the homotopy generators //(gamma_1,...,gamma_r)=(\alpha_1,...,\alpha_{p+q},\alpha,\beta))
local B;
//i=+/-1,...,=/-#W-1
B:=W;
if j ge 1 then
B[j]:=W[j+1];
B[j+1]:=Conjugate(W[j+1],W[j]);
else 
B[-j]:=ConjugateInverse(W[-j],W[-j+1]);
B[-j+1]:=W[-j];
end if;
return B;
end function;

BraidAction:=function(W,Word)
//W=[[[1,1]],...,[[r,1]]] Tuple describing the homotopy generators 
//(gamma_1,...,gamma_r)=(\alpha_1,...,\alpha_{p+q},\alpha,\beta)
//Word of braids  in terms of a list [i_1,....,i_k] where 
// a positive index k  stands for beta_k (cf. Section 2.1) and
// a negative index -k (k>0) stands for \beta_k^-1
local BA,i;
BA:=W;
for i in [1..#(Word)] do
BA:=Braid(Word[i],BA);
end for;	
return BA;
end function;

MonodromyOfDeformations:=function(W,A,K)
//K field
//A\in \GL_{n_1n_2}(K)^{p+q+2} monodromy tuple of L_1\otimes L_2(y_0-x)
// as in Lemma 3.3.1
// W=[W[1],...,W[r]] List of words in the homotopy generators 
// (gamma_1,...,gamma_r)=(\alpha_1,...,\alpha_{p+q},\alpha,\beta)
//which describe the effect of deformation along  a list of braids, e.g., obtained by 
// an output of BraidAction(W0,delta_{i,j})
//Output:  a matrix  B in GL(K^nr) whose restriction to the the space H_A, described in 
//Lemma 3.2.1, 
//gives the map on cocycles as in Thm. 3.3.3., induced by 
// a deformation of  (gamma_1,...,gamma_r) given by inverse conjugation of some braid,
// described by W, e.g.,  \delta\mapsto \delta^{delta_{i,j}}
local r,V,i,j,M,Matrixx,k,n;
n:=Rank(A[1]);
r:=#(A);
V:=[];
for i in [1..r] do 
V[i]:=[]; V1:=[];
V[i][#(W[i])]:=IdentityMatrix(K,n);
for j in [1..#W[i]-1] do 
V[i][#(W[i])-j]:=A[W[i][#W[i]+1-j][1]]^W[i][#W[i]+1-j][2]*V[i][#W[i]+1-j];
end for;
end for;

//V[i]=[[W[i],A[W[i][2][1]]*...*A[W[i][#[W[i]]][1]], ... , [W[i][#[i]],IdentityMatrix]]
for i in [1..r] do
for j in [1..#(V[i])] do
if W[i][j][2] eq -1 then 
V[i][j]:=-A[W[i][j][1]]^-1*V[i][j];
end if;
end for;
end for;

Matrixx:=ZeroMatrix(K,n*r,n*r);
for i in [1..r] do
for j in [1..r] do
M:=ZeroMatrix(K,n,n);
for k in [1..#V[i]] do 
if W[i][k][1] eq j then
M:=M+V[i][k];
end if;
end for;
Matrixx:=InsertBlock(Matrixx,M,(j-1)*n+1,(i-1)*n+1);
end for;                 
end for;                     
return Matrixx;
end function;

// In the following, we construct the G2-example associated to a 7-sheaf:
///////////////////////////////////////////////////////////////////////////////////////
//the four cases in  Tab. 1
K<y>:=FunctionField(Rationals());
B:=[];
B[1]:=Matrix(K,[ [ 1, 1 ], [ 0, 1 ] ]); //#0
B[2]:=Matrix(K,[ [ 1, -1 ], [ 0, 1 ] ]);  //#1
B[3]:=Matrix(K,[ [ 1, 0 ], [ y, 1 ] ]); //#t
B[4]:=(B[1]*B[2]*B[3])^-1; //#infty

K<a,b>:=FunctionField(Rationals(),2);
B:=[];
B[1]:=Matrix(K,[ [ 1, 1 ], [ 0, 1 ] ]); //#0
B[2]:=Matrix(K,[ [ 1, 0 ], [ (a-1)^2/(b*(a+b)), 1 ] ]);  //#1
B[3]:=Matrix(K,[ [ a, b ], [ -(a-1)^2/b, 2-a ] ]); //#t
B[4]:=(B[1]*B[2]*B[3])^-1; //#infty

K<y>:=FunctionField(Rationals());
B:=[];
B[1]:=Matrix(K,[ [ 1, 1 ], [ 0, 1 ] ]); //#0
B[2]:=Matrix(K,[ [ 1, 0 ], [ y, 1 ] ]);  //#1
B[3]:=Matrix(K,[ [ 1, 0 ], [ -y, 1 ] ]); //#t
B[4]:=(B[1]*B[2]*B[3])^-1; //#infty

K<y>:=FunctionField(Rationals());
B:=[];
B[1]:=Matrix(K,[ [ 1, 1 ], [ 0, 1 ] ]); //#0
B[2]:=Matrix(K,[ [ 1, 0 ], [ y, 1 ] ]);  //#1
B[3]:=Matrix(K,[ [ 1, -1 ], [ 0, 1 ] ]); //#t
B[4]:=(B[1]*B[2]*B[3])^-1; //#infty

//////////////////////////////////////////////////////////////////////////////////////////////
A:=[B[4],B[4]^-1*B[3]*B[4],(B[3]*B[4])^-1*B[2]*(B[3]*B[4]),B[1]^2,B[2],B[3],B[4]];
//A is the singular monodromy tuple of the 7-point sheaf of Eq. 4.1
AT:=[];
for i in [1..7] do
AT[i]:=KroneckerProduct(A[i],I2);
end for;

P:=[];
P[1]:=A[7];
for i in [1 .. 6]  do
P[i+1]:=A[7-i]*P[i];
end for;
for i in [1 .. 7] do
AT[i+7]:=KroneckerProduct(I2,P[i]^-1*A[7-i+1]*P[i]);
end for;
AT[15]:=AT[1]^0;
AT[16]:=AT[1]^0;
//AT is the monodromy tuple of N\otimes N(y_0-x) 
//WW0 is represents the tuple of  homotopy generators 
// (gamma_1,...,gamma_16)=(\alpha_1,...,\alpha_{7+7},\alpha,\beta))
WW0:=[[[1,1]],[[2,1]],[[3,1]],[[4,1]],[[5,1]],[[6,1]],[[7,1]],[[8,1]],[[9,1]], [[10,1]],[[11,1]],[[12,1]],[[13,1]],[[14,1]],[[15,1]],[[16,1]]];
WW:=[];
//WW19--WW1 gives the deformation along \delta_{19},...,\delta_1 as in Prop 2.3.1
//and Remark 2.3.2.
for i in [1.. 19] do
WW[i]:=BraidAction(WW0,delta[i]);
end for;

//WW20 gives the deformation along \hat{\alpha} as in Prop2.3.3 
WW[20]:=[[[14,-1],[13,-1],[12,-1],[11,-1],[10,-1],[9,-1],[8,-1],[1,1],[8,1],[9,1], [10,1],[11,1],[12,1],[13,1],[14,1]],[[14,-1],[13,-1],[12,-1],[11,-1],[10,-1], [9,-1],[8,-1],[2,1],[8,1],[9,1],[10,1],[11,1],[12,1],[13,1],[14,1]],
[[14,-1],[13,-1],[12,-1],[11,-1],[10,-1],[9,-1],[8,-1],[3,1],[8,1],[9,1],[10,1], [11,1],[12,1],[13,1],[14,1]],
[[14,-1],[13,-1],[12,-1],[11,-1],[10,-1],[9,-1],[8,-1],[4,1],[8,1],[9,1],[10,1], [11,1],[12,1],[13,1],[14,1]],
[[14,-1],[13,-1],[12,-1],[11,-1],[10,-1],[9,-1],[8,-1],[5,1],[8,1],[9,1],[10,1], [11,1],[12,1],[13,1],[14,1]],
[[14,-1],[13,-1],[12,-1],[11,-1],[10,-1],[9,-1],[8,-1],[6,1],[8,1],[9,1],[10,1], [11,1],[12,1],[13,1],[14,1]],
[[14,-1],[13,-1],[12,-1],[11,-1],[10,-1],[9,-1],[8,-1],[7,1],[8,1],[9,1],[10,1],
[11,1],[12,1],[13,1],[14,1]],
[[15,1],[8,1],[15,-1]], [[15,1],[9,1],[15,-1]], [[15,1],[10,1],[15,-1]],
[[15,1],[11,1],[15,-1]], [[15,1],[12,1],[15,-1]], [[15,1],[13,1],[15,-1]],
[[15,1],[14,1],[15,-1]], [[15,1]], 
[[14,-1],[13,-1],[12,-1],[11,-1],[10,-1],[9,-1],[8,-1],[16,1]]];
//WW21 gives the deformation along \hat{\beta} as in Prop 2.3.3
WW[21]:=[[[16,1],[8,1],[9,1],[10,1],[11,1],[12,1],[13,1],[14,1],[16,-1],[1,1], [16,1],[14,-1],[13,-1],[12,-1],[11,-1],[10,-1],[9,-1],[8,-1],[16,-1]],
[[16,1],[8,1],[9,1],[10,1],[11,1],[12,1],[13,1],[14,1],[16,-1],[2,1],[16,1], [14,-1],[13,-1],[12,-1],[11,-1],[10,-1],[9,-1],[8,-1],[16,-1]],
[[16,1],[8,1],[9,1],[10,1],[11,1],[12,1],[13,1],[14,1],[16,-1],[3,1],[16,1], [14,-1],[13,-1],[12,-1],[11,-1],[10,-1],[9,-1],[8,-1],[16,-1]],
[[16,1],[8,1],[9,1],[10,1],[11,1],[12,1],[13,1],[14,1],[16,-1],[4,1],[16,1], [14,-1],[13,-1],[12,-1],[11,-1],[10,-1],[9,-1],[8,-1],[16,-1]],
[[16,1],[8,1],[9,1],[10,1],[11,1],[12,1],[13,1],[14,1],[16,-1],[5,1],[16,1], [14,-1],[13,-1],[12,-1],[11,-1],[10,-1],[9,-1],[8,-1],[16,-1]],
[[16,1],[8,1],[9,1],[10,1],[11,1],[12,1],[13,1],[14,1],[16,-1],[6,1],[16,1], [14,-1],[13,-1],[12,-1],[11,-1],[10,-1],[9,-1],[8,-1],[16,-1]],
[[16,1],[8,1],[9,1],[10,1],[11,1],[12,1],[13,1],[14,1],[16,-1],[7,1],[16,1], [14,-1],[13,-1],[12,-1],[11,-1],[10,-1],[9,-1],[8,-1],[16,-1]],
[[16,1],[8,1],[16,-1]], [[16,1],[9,1],[16,-1]], [[16,1],[10,1],[16,-1]],
[[16,1],[11,1],[16,-1]], [[16,1],[12,1],[16,-1]], [[16,1],[13,1],[16,-1]],
[[16,1],[14,1],[16,-1]], [[8,1],[9,1],[10,1],[11,1],[12,1],[13,1],[14,1],[15,1]],
[[16,1]]];

//MM[1] to MM[21] give the Monodromy of \delta_{1},...\delta_{19}, \hat{\alpha}^{-1}, 
//\hat{\beta}^{-1} (using that the global
//monodromy of AT is trivial 
// Note that one has to invert MM[20] and MM[21] in order to obtain the elliptic product relation

MM:=[];
for i in [1.. 21] do
MM[i]:=MonodromyOfDeformations(WW[i],AT,K);
end for;

//To analyze further, one may proceed as follows:

MatAlg:=MatrixAlgebra<K,Rank(MM[1])|MM>;
TT:=RModule(MatAlg);
C:=ConstituentsWithMultiplicities(TT);
C;

//A final remark: When dealing with all evaluations, then the spaces 
//ker(A_i-1)=\mathbb{Q} contribute a shift of  L 
// because one is actually computing the variation of H^1(E,j_!L\otimes L(y_0-x)) using 
//the fact that locally at x_i the extension by zero 
//j_!L is an extension of j_*L by \delta_{x_i}, explaining the rest of the occurring 
//rank-two factors in C. 
\end{lstlisting}

\section{Appendix -- Monodromy for seven-point sheaves of Beauville type}\label{app:BeauvMon}

We give the explicit monodromy matrices of a seven-point sheaf $N$ on $E$ obtained as in \cite{KatzEll}, \S~4 by the Beauville classification. We refer to ibid. and \cite{CDR18ARXIV},\S~4.1 for details.

\bigskip

In the case of the Beauville family $f\colon\mathcal{B}\to \PP^1_\lambda$ defined by $y^2=x(x-1)(\lambda^2-x)$, the monodromy tuple of $\mathcal G:=R^1f_*(\bQl)$ at $\mathbf{x}=\{-1,0,1,\infty\}\in \PP^1$ of \cite{KatzEll} is given by:
\begin{equation}
\Bigl(M_{-1},M_1,M_0,M_{\infty}\Bigr)=
\left(\begin{pmatrix}
1& 0\\
2&1
\end{pmatrix},
\begin{pmatrix}
-19& -8\\
50&21
\end{pmatrix},
\begin{pmatrix}
-7& -4\\
16&9
\end{pmatrix},
\begin{pmatrix}
-3& -4\\
4&5
\end{pmatrix}\right).
\end{equation}

\bigskip

On the generic elliptic curve $E=E_t\subset B$ defined by $y^2=t(x^3-x)+t^2$, the pullback $L$ of $\mathcal G$ along the first projection $E\to \PP^1_\lambda$ gives a seven-point sheaf $N$ on $E$ as in Lem.~\ref{lem:7sheafP1} whose local monodromy matrices at $\mathbf{a}=\{-c,-b,-a,0,a,b,c\}\in E$ is given  by:
\begin{multline}
\Bigl(M_{-c},M_{-b},M_{-a},M_{0},M_{a},M_{b},M_{c}\Bigr)=\\
\left(\left(\begin{array}{cc}
-3& -4\\
4&5
\end{array}
\right), 
\left(\begin{array}{cc}
-23& -36\\
16&25
\end{array}
\right), 
\left(\begin{array}{cc}
-3& -8\\
2&5
\end{array}
\right), 
\left(\begin{array}{cc}
1& 0\\
4&1
\end{array}
\right),\right.\\ 
\left.\left(\begin{array}{cc}
-19& -8\\
50&21
\end{array}
\right),
\left(\begin{array}{cc}
-7& -4\\
16&9
\end{array}
\right), 
\left(\begin{array}{cc}
-3& -4\\
4&5
\end{array}
\right)\right). 
\end{multline}

\pagestyle{plain}


\newcommand{\SortNoop}[1]{}\def\cprime{$'$}
\providecommand{\bysame}{\leavevmode\hbox to3em{\hrulefill}\thinspace}
\providecommand{\MR}{\relax\ifhmode\unskip\space\fi MR }
\providecommand{\MRhref}[2]{%
	\href{http://www.ams.org/mathscinet-getitem?mr=#1}{#2}
}
\providecommand{\href}[2]{#2}

\end{document}

%% file: Figures/Braids_Conf.pdf_tex
\begingroup%
  \makeatletter%
  \providecommand\color[2][]{%
    \errmessage{(Inkscape) Color is used for the text in Inkscape, but the package 'color.sty' is not loaded}%
    \renewcommand\color[2][]{}%
  }%
  \providecommand\transparent[1]{%
    \errmessage{(Inkscape) Transparency is used (non-zero) for the text in Inkscape, but the package 'transparent.sty' is not loaded}%
    \renewcommand\transparent[1]{}%
  }%
  \providecommand\rotatebox[2]{#2}%
  \ifx\svgwidth\undefined%
    \setlength{\unitlength}{432.34021284bp}%
    \ifx\svgscale\undefined%
      \relax%
    \else%
      \setlength{\unitlength}{\unitlength * \real{\svgscale}}%
    \fi%
  \else%
    \setlength{\unitlength}{\svgwidth}%
  \fi%
  \global\let\svgwidth\undefined%
  \global\let\svgscale\undefined%
  \makeatother%
  \begin{picture}(1,1)%
    \put(0,0){\includegraphics[width=\unitlength,page=1]{Braids_Conf.pdf}}%
    \put(0.07610314,0.78739749){\color[rgb]{0,0,0}\makebox(0,0)[lb]{\smash{$x_1$}}}%
    \put(0.2337171,0.63151836){\color[rgb]{0,0,0}\makebox(0,0)[lb]{\smash{$x_p$}}}%
    \put(0.26940327,0.55320692){\color[rgb]{0,0,0}\makebox(0,0)[lb]{\smash{$x_{p+1}$}}}%
    \put(0.48049516,0.33854074){\color[rgb]{0,0,0}\makebox(0,0)[lb]{\smash{$x_{p+q}$}}}%
    \put(0.74938072,0.16382733){\color[rgb]{0,0,0}\makebox(0,0)[lb]{\smash{$x_{0}$}}}%
    \put(0.22094097,0.91854168){\color[rgb]{0,0,0}\makebox(0,0)[lb]{\smash{$\alpha_1$}}}%
    \put(0.62285579,0.28608651){\color[rgb]{0,0,0}\makebox(0,0)[lb]{\smash{$\alpha_{p+q}$}}}%
  \end{picture}%
\endgroup%

%% file: Figures/Braids_Act1.pdf_tex
\begingroup%
  \makeatletter%
  \providecommand\color[2][]{%
    \errmessage{(Inkscape) Color is used for the text in Inkscape, but the package 'color.sty' is not loaded}%
    \renewcommand\color[2][]{}%
  }%
  \providecommand\transparent[1]{%
    \errmessage{(Inkscape) Transparency is used (non-zero) for the text in Inkscape, but the package 'transparent.sty' is not loaded}%
    \renewcommand\transparent[1]{}%
  }%
  \providecommand\rotatebox[2]{#2}%
  \ifx\svgwidth\undefined%
    \setlength{\unitlength}{432.34021284bp}%
    \ifx\svgscale\undefined%
      \relax%
    \else%
      \setlength{\unitlength}{\unitlength * \real{\svgscale}}%
    \fi%
  \else%
    \setlength{\unitlength}{\svgwidth}%
  \fi%
  \global\let\svgwidth\undefined%
  \global\let\svgscale\undefined%
  \makeatother%
  \begin{picture}(1,1.00312409)%
    \put(0,0){\includegraphics[width=\unitlength,page=1]{Braids_Act1.pdf}}%
    \put(0.25205586,0.65853071){\color[rgb]{0,0,0}\makebox(0,0)[lb]{\smash{$x_p$}}}%
    \put(0.1280082,0.76613188){\color[rgb]{0,0,0}\makebox(0,0)[lb]{\smash{$x_1$}}}%
    \put(0.34556865,0.56863933){\color[rgb]{0,0,0}\makebox(0,0)[lb]{\smash{$x_{p+1}$}}}%
    \put(0.54123373,0.35415333){\color[rgb]{0,0,0}\makebox(0,0)[lb]{\smash{$x_{p+q}$}}}%
    \put(0.74196852,0.17943972){\color[rgb]{0,0,0}\makebox(0,0)[lb]{\smash{$x_{0}$}}}%
  \end{picture}%
\endgroup%

%% file: Figures/Braids_Act2.pdf_tex
\begingroup%
  \makeatletter%
  \providecommand\color[2][]{%
    \errmessage{(Inkscape) Color is used for the text in Inkscape, but the package 'color.sty' is not loaded}%
    \renewcommand\color[2][]{}%
  }%
  \providecommand\transparent[1]{%
    \errmessage{(Inkscape) Transparency is used (non-zero) for the text in Inkscape, but the package 'transparent.sty' is not loaded}%
    \renewcommand\transparent[1]{}%
  }%
  \providecommand\rotatebox[2]{#2}%
  \ifx\svgwidth\undefined%
    \setlength{\unitlength}{432.34021284bp}%
    \ifx\svgscale\undefined%
      \relax%
    \else%
      \setlength{\unitlength}{\unitlength * \real{\svgscale}}%
    \fi%
  \else%
    \setlength{\unitlength}{\svgwidth}%
  \fi%
  \global\let\svgwidth\undefined%
  \global\let\svgscale\undefined%
  \makeatother%
  \begin{picture}(1,1.00312409)%
    \put(0,0){\includegraphics[width=\unitlength,page=1]{Braids_Act2.pdf}}%
    \put(0.75790071,0.16934682){\color[rgb]{0,0,0}\makebox(0,0)[lb]{\smash{$x_{0}$}}}%
    \put(0.57032792,0.34653524){\color[rgb]{0,0,0}\makebox(0,0)[lb]{\smash{$x_{p+q}$}}}%
    \put(0.11306858,0.79274336){\color[rgb]{0,0,0}\makebox(0,0)[lb]{\smash{$x_1$}}}%
    \put(0.33413251,0.54630389){\color[rgb]{0,0,0}\makebox(0,0)[lb]{\smash{$x_{p+1}$}}}%
    \put(0.27376366,0.62002409){\color[rgb]{0,0,0}\makebox(0,0)[lb]{\smash{$x_p$}}}%
  \end{picture}%
\endgroup%

%% file: Figures/Braids_Conf2.pdf_tex
\begingroup%
  \makeatletter%
  \providecommand\color[2][]{%
    \errmessage{(Inkscape) Color is used for the text in Inkscape, but the package 'color.sty' is not loaded}%
    \renewcommand\color[2][]{}%
  }%
  \providecommand\transparent[1]{%
    \errmessage{(Inkscape) Transparency is used (non-zero) for the text in Inkscape, but the package 'transparent.sty' is not loaded}%
    \renewcommand\transparent[1]{}%
  }%
  \providecommand\rotatebox[2]{#2}%
  \ifx\svgwidth\undefined%
    \setlength{\unitlength}{432.34021284bp}%
    \ifx\svgscale\undefined%
      \relax%
    \else%
      \setlength{\unitlength}{\unitlength * \real{\svgscale}}%
    \fi%
  \else%
    \setlength{\unitlength}{\svgwidth}%
  \fi%
  \global\let\svgwidth\undefined%
  \global\let\svgscale\undefined%
  \makeatother%
  \begin{picture}(1,1)%
    \put(0,0){\includegraphics[width=\unitlength,page=1]{Braids_Conf2.pdf}}%
    \put(0.0277781,0.80598411){\color[rgb]{0,0,0}\makebox(0,0)[lb]{\smash{$x_1$}}}%
    \put(0.39499704,0.44014714){\color[rgb]{0,0,0}\makebox(0,0)[lb]{\smash{$x_{r-1}$}}}%
    \put(0.82992234,0.14152341){\color[rgb]{0,0,0}\makebox(0,0)[lb]{\smash{$x_{0}$}}}%
    \put(0.39689368,0.8342826){\color[rgb]{0,0,0}\makebox(0,0)[lb]{\smash{$\alpha_1$}}}%
    \put(0.60922567,0.30467308){\color[rgb]{0,0,0}\makebox(0,0)[lb]{\smash{$\alpha_{r}$}}}%
    \put(0.08384579,0.17580636){\color[rgb]{0,0,0}\makebox(0,0)[lb]{\smash{$\alpha$}}}%
    \put(0,0){\includegraphics[width=\unitlength,page=2]{Braids_Conf2.pdf}}%
    \put(0.1064256,0.76184788){\color[rgb]{0,0,0}\makebox(0,0)[lb]{\smash{$x_2$}}}%
    \put(0.46593194,0.38188562){\color[rgb]{0,0,0}\makebox(0,0)[lb]{\smash{$x_{r}$}}}%
    \put(0,0){\includegraphics[width=\unitlength,page=3]{Braids_Conf2.pdf}}%
    \put(0.81603967,0.80842149){\color[rgb]{0,0,0}\makebox(0,0)[lb]{\smash{$\beta$}}}%
  \end{picture}%
\endgroup%

%% file: Figures/Braids_Conf_Base.pdf_tex
\begingroup%
  \makeatletter%
  \providecommand\color[2][]{%
    \errmessage{(Inkscape) Color is used for the text in Inkscape, but the package 'color.sty' is not loaded}%
    \renewcommand\color[2][]{}%
  }%
  \providecommand\transparent[1]{%
    \errmessage{(Inkscape) Transparency is used (non-zero) for the text in Inkscape, but the package 'transparent.sty' is not loaded}%
    \renewcommand\transparent[1]{}%
  }%
  \providecommand\rotatebox[2]{#2}%
  \ifx\svgwidth\undefined%
    \setlength{\unitlength}{432.34021284bp}%
    \ifx\svgscale\undefined%
      \relax%
    \else%
      \setlength{\unitlength}{\unitlength * \real{\svgscale}}%
    \fi%
  \else%
    \setlength{\unitlength}{\svgwidth}%
  \fi%
  \global\let\svgwidth\undefined%
  \global\let\svgscale\undefined%
  \makeatother%
  \begin{picture}(1,1)%
    \put(0,0){\includegraphics[width=\unitlength,page=1]{Braids_Conf_Base.pdf}}%
    \put(0.02282169,0.75889817){\color[rgb]{0,0,0}\makebox(0,0)[lb]{\smash{$x_1+y_q$}}}%
    \put(0.82992234,0.14152341){\color[rgb]{0,0,0}\makebox(0,0)[lb]{\smash{$y_{0}$}}}%
    \put(0.39689368,0.8342826){\color[rgb]{0,0,0}\makebox(0,0)[lb]{\smash{$\delta_{1,q}$}}}%
    \put(0.80624311,0.79164078){\color[rgb]{0,0,0}\makebox(0,0)[lb]{\smash{$\hat{\beta}$}}}%
    \put(0,0){\includegraphics[width=\unitlength,page=2]{Braids_Conf_Base.pdf}}%
    \put(0.45849733,0.31745236){\color[rgb]{0,0,0}\makebox(0,0)[lb]{\smash{$x_p+y_1$}}}%
    \put(0,0){\includegraphics[width=\unitlength,page=3]{Braids_Conf_Base.pdf}}%
    \put(0.13948914,0.14426203){\color[rgb]{0,0,0}\makebox(0,0)[lb]{\smash{$\hat{\alpha}$}}}%
    \put(0.62733068,0.2680885){\color[rgb]{0,0,0}\makebox(0,0)[lb]{\smash{$\delta_{p,1}$}}}%
  \end{picture}%
\endgroup%

%% file: Figures/Phi_Base.pdf_tex
\begingroup%
  \makeatletter%
  \providecommand\color[2][]{%
    \errmessage{(Inkscape) Color is used for the text in Inkscape, but the package 'color.sty' is not loaded}%
    \renewcommand\color[2][]{}%
  }%
  \providecommand\transparent[1]{%
    \errmessage{(Inkscape) Transparency is used (non-zero) for the text in Inkscape, but the package 'transparent.sty' is not loaded}%
    \renewcommand\transparent[1]{}%
  }%
  \providecommand\rotatebox[2]{#2}%
  \ifx\svgwidth\undefined%
    \setlength{\unitlength}{432.34021284bp}%
    \ifx\svgscale\undefined%
      \relax%
    \else%
      \setlength{\unitlength}{\unitlength * \real{\svgscale}}%
    \fi%
  \else%
    \setlength{\unitlength}{\svgwidth}%
  \fi%
  \global\let\svgwidth\undefined%
  \global\let\svgscale\undefined%
  \makeatother%
  \begin{picture}(1,1)%
    \put(0,0){\includegraphics[width=\unitlength,page=1]{Phi_Base.pdf}}%
    \put(0.23763084,0.54423203){\color[rgb]{0,0,0}\makebox(0,0)[lb]{\smash{$\phi(x_0)=-x_{0}$}}}%
    \put(0.46752259,0.08834215){\color[rgb]{0,0,0}\makebox(0,0)[lb]{\smash{$\phi(\gamma_{1})$}}}%
    \put(0.23873373,0.87837806){\color[rgb]{0,0,0}\makebox(0,0)[lb]{\smash{$\phi(\beta)$}}}%
    \put(0,0){\includegraphics[width=\unitlength,page=2]{Phi_Base.pdf}}%
    \put(0.49319201,0.42153716){\color[rgb]{0,0,0}\makebox(0,0)[lb]{\smash{$-x_r$}}}%
    \put(0,0){\includegraphics[width=\unitlength,page=3]{Phi_Base.pdf}}%
    \put(0.71691137,0.52218853){\color[rgb]{0,0,0}\makebox(0,0)[lb]{\smash{$\phi(\alpha)$}}}%
    \put(0.867384,0.23646794){\color[rgb]{0,0,0}\makebox(0,0)[lb]{\smash{$-x_1$}}}%
    \put(0.2970005,0.41804791){\color[rgb]{0,0,0}\makebox(0,0)[lb]{\smash{$\phi(\gamma_{r})$}}}%
  \end{picture}%
\endgroup%

%% file: Figures/Phi_Base2.pdf_tex
\begingroup%
  \makeatletter%
  \providecommand\color[2][]{%
    \errmessage{(Inkscape) Color is used for the text in Inkscape, but the package 'color.sty' is not loaded}%
    \renewcommand\color[2][]{}%
  }%
  \providecommand\transparent[1]{%
    \errmessage{(Inkscape) Transparency is used (non-zero) for the text in Inkscape, but the package 'transparent.sty' is not loaded}%
    \renewcommand\transparent[1]{}%
  }%
  \providecommand\rotatebox[2]{#2}%
  \ifx\svgwidth\undefined%
    \setlength{\unitlength}{432.34021284bp}%
    \ifx\svgscale\undefined%
      \relax%
    \else%
      \setlength{\unitlength}{\unitlength * \real{\svgscale}}%
    \fi%
  \else%
    \setlength{\unitlength}{\svgwidth}%
  \fi%
  \global\let\svgwidth\undefined%
  \global\let\svgscale\undefined%
  \makeatother%
  \begin{picture}(1,1)%
    \put(0,0){\includegraphics[width=\unitlength,page=1]{Phi_Base2.pdf}}%
    \put(0.38809472,0.41168658){\color[rgb]{0,0,0}\makebox(0,0)[lb]{\smash{$x_1'$}}}%
    \put(0,0){\includegraphics[width=\unitlength,page=2]{Phi_Base2.pdf}}%
    \put(0.29957889,0.06524048){\color[rgb]{0,0,0}\makebox(0,0)[lb]{\smash{$\alpha'$}}}%
    \put(0.70616719,0.24435355){\color[rgb]{0,0,0}\makebox(0,0)[lb]{\smash{$x_r'$}}}%
    \put(0,0){\includegraphics[width=\unitlength,page=3]{Phi_Base2.pdf}}%
    \put(0.44627744,0.59327275){\color[rgb]{0,0,0}\makebox(0,0)[lb]{\smash{$\gamma_1'$}}}%
    \put(0.56688359,0.35571618){\color[rgb]{0,0,0}\makebox(0,0)[lb]{\smash{$\gamma_r'$}}}%
    \put(0,0){\includegraphics[width=\unitlength,page=4]{Phi_Base2.pdf}}%
    \put(0.8878222,0.87400077){\color[rgb]{0,0,0}\makebox(0,0)[lb]{\smash{$\beta'$}}}%
  \end{picture}%
\endgroup%

%% file: Figures/Convol.pdf_tex
\begingroup%
  \makeatletter%
  \providecommand\color[2][]{%
    \errmessage{(Inkscape) Color is used for the text in Inkscape, but the package 'color.sty' is not loaded}%
    \renewcommand\color[2][]{}%
  }%
  \providecommand\transparent[1]{%
    \errmessage{(Inkscape) Transparency is used (non-zero) for the text in Inkscape, but the package 'transparent.sty' is not loaded}%
    \renewcommand\transparent[1]{}%
  }%
  \providecommand\rotatebox[2]{#2}%
  \ifx\svgwidth\undefined%
    \setlength{\unitlength}{525.60068686bp}%
    \ifx\svgscale\undefined%
      \relax%
    \else%
      \setlength{\unitlength}{\unitlength * \real{\svgscale}}%
    \fi%
  \else%
    \setlength{\unitlength}{\svgwidth}%
  \fi%
  \global\let\svgwidth\undefined%
  \global\let\svgscale\undefined%
  \makeatother%
  \begin{picture}(1,0.82256402)%
    \put(0,0){\includegraphics[width=\unitlength,page=1]{Convol.pdf}}%
    \put(0.3282725,0.0289978){\color[rgb]{0,0,0}\makebox(0,0)[lb]{\smash{$x_1$}}}%
    \put(0,0){\includegraphics[width=\unitlength,page=2]{Convol.pdf}}%
    \put(0.37821537,0.0289978){\color[rgb]{0,0,0}\makebox(0,0)[lb]{\smash{$x_2$}}}%
    \put(0.45738873,0.02830604){\color[rgb]{0,0,0}\makebox(0,0)[lb]{\smash{$x_3$}}}%
    \put(0.5094937,0.0289978){\color[rgb]{0,0,0}\makebox(0,0)[lb]{\smash{$x_4$}}}%
    \put(0.55943654,0.0289978){\color[rgb]{0,0,0}\makebox(0,0)[lb]{\smash{$x_5$}}}%
    \put(0.63506443,0.027571){\color[rgb]{0,0,0}\makebox(0,0)[lb]{\smash{$x_6$}}}%
    \put(0.68643444,0.027571){\color[rgb]{0,0,0}\makebox(0,0)[lb]{\smash{$x_7$}}}%
    \put(0,0){\includegraphics[width=\unitlength,page=3]{Convol.pdf}}%
    \put(0.54725689,0.10980105){\color[rgb]{0,0,0}\makebox(0,0)[lb]{\smash{$pr_1$}}}%
    \put(0,0){\includegraphics[width=\unitlength,page=4]{Convol.pdf}}%
    \put(0.81490197,0.43228907){\color[rgb]{0,0,0}\makebox(0,0)[lb]{\smash{$pr_2$}}}%
    \put(0.18248918,0.609848){\color[rgb]{0,0,0}\rotatebox{30}{\makebox(0,0)[lb]{\smash{$\mathbb{V}\subset E_x\times E_x$}}}}%
    \put(0.82863782,0.65862146){\color[rgb]{0,0,0}\makebox(0,0)[lb]{\smash{$E\setminus \x_1\ast \x_2$}}}%
    \put(0.80961069,0.05120934){\color[rgb]{0,0,0}\makebox(0,0)[lb]{\smash{$E\setminus \x_1$}}}%
    \put(0,0){\includegraphics[width=\unitlength,page=5]{Convol.pdf}}%
    \put(0.91581476,0.60566851){\color[rgb]{0,0,0}\makebox(0,0)[lb]{\smash{$x_7+y_1$}}}%
    \put(0.91587845,0.57854712){\color[rgb]{0,0,0}\makebox(0,0)[lb]{\smash{$x_6+y_1$}}}%
    \put(0,0){\includegraphics[width=\unitlength,page=6]{Convol.pdf}}%
    \put(0.18008429,0.23543927){\color[rgb]{0,0,0}\makebox(0,0)[lb]{\smash{$d$}}}%
    \put(0,0){\includegraphics[width=\unitlength,page=7]{Convol.pdf}}%
    \put(0.0179509,0.299382){\color[rgb]{0,0,0}\makebox(0,0)[lb]{\smash{$E\setminus \x_2$}}}%
    \put(0.0966146,0.07625596){\color[rgb]{0,0,0}\makebox(0,0)[lb]{\smash{$y_7$}}}%
    \put(0.08477122,0.09892126){\color[rgb]{0,0,0}\makebox(0,0)[lb]{\smash{$y_6$}}}%
    \put(0.02217267,0.21309063){\color[rgb]{0,0,0}\makebox(0,0)[lb]{\smash{$y_2$}}}%
    \put(0.01137201,0.23536644){\color[rgb]{0,0,0}\makebox(0,0)[lb]{\smash{$y_1$}}}%
    \put(0,0){\includegraphics[width=\unitlength,page=8]{Convol.pdf}}%
    \put(0.0407226,0.17385373){\color[rgb]{0,0,0}\makebox(0,0)[lb]{\smash{$y_3$}}}%
    \put(0.05572717,0.15456845){\color[rgb]{0,0,0}\makebox(0,0)[lb]{\smash{$y_4$}}}%
    \put(0.06714266,0.13316437){\color[rgb]{0,0,0}\makebox(0,0)[lb]{\smash{$y_5$}}}%
    \put(0.91835662,0.19894602){\color[rgb]{0,0,0}\makebox(0,0)[lb]{\smash{$x_1+y_7$}}}%
    \put(0.91927073,0.2220223){\color[rgb]{0,0,0}\makebox(0,0)[lb]{\smash{$x_1+y_6$}}}%
  \end{picture}%
\endgroup%